\documentclass[a4,12pt]{article}
\usepackage{amsmath}
\usepackage{amssymb}
\usepackage{graphicx}
\usepackage{float}
\usepackage{cite}
\usepackage{subcaption}
\usepackage[T1]{fontenc}
\usepackage{pgfplotstable}
\usepackage{booktabs}
\usepackage[top=1in, bottom=0.5in, left=1in, right=1in]{geometry}
\begin{document}
\begin{center}
\textbf{\Large{Determination of forcing functions in the wave equation. Part I: the space-dependent case}}
\end{center}
S.O. Hussein and D. Lesnic 
\\
\textit{Department of Applied Mathematics, University of Leeds, Leeds LS2 9JT, UK}\\
E-mails: ml10soh@leeds.ac.uk, D.Lesnic@leeds.ac.uk\\
\\
\textbf{\large Abstract.}
We consider the inverse problem for the wave equation which consists of determining an unknown space-dependent force function acting on a vibrating structure from Cauchy boundary data. Since only boundary data are used as measurements, the study has importance and significance to non-intrusive and non-destructive testing of materials. This inverse force problem is linear, the solution is unique, but the problem is still ill-posed since, in general, the solution does not exist and, even if it exists, it does not depend continuously upon the input data. Numerically, the finite difference method combined with the Tikhonov regularization are employed in order to obtain a stable solution. Several orders of regularization are investigated. The choice of the regularization parameter is based on the L-curve method. Numerical results show that the solution is accurate for exact data and stable for noisy data. An extension to the case of multiple additive forces is also addressed. In a companion paper, in Part II, the time-dependent force identification will be undertaken.

\section{Introduction}
The aim of this paper is to investigate an inverse force problem for the hyperbolic wave equation. The forcing function is assumed to depend only upon the space variable in order to ensure uniqueness of the solution, \cite{candunn70,engl94,vmk92,my95}. These authors have given conditions to be satisfied by the data in order to ensure uniqueness and, in the case of \cite{candunn70}, continuous dependence upon the data. However, no numerical results were presented and it is the main purpose and novelty of our study to develop an efficient numerical solution for this inverse linear, but ill-posed problem. In a previous study, \cite{hussein14}, we have used the boundary element method (BEM) to numerically discretise the wave equation with constant wave speed based on the available fundamental solution, \cite{morfesh53}. Furthermore, by assuming that the force function $f(x)$ appears as a free term in the wave equation, the method of separating variables, \cite{candunn70}, was applicable and regularisation was used to stabilise the resulting system of linear algebraic equations. However, if the wave speed is not constant or, if the force appears in a non-free term as $f(x)h(x,t)$ the above methods are not applicable. Therefore, in order to extend this range of applicability, in this paper the numerical method for discretising the wave equation is the finite difference method (FDM). The resulting system of linear equations is ill-conditioned, the original problem being ill-posed. The choice of the regularization parameter introduced by this technique is important for the stability of the numerical solution and in our study this is based on the L-curve criterion,\cite{hansen2001}. 

The structure of the paper is as follows. In Section 2, we briefly describe inverse force problems for the hyperbolic wave equation recalling the uniqueness theorems of \cite{engl94,vmk92,my95}. In Sections 3 and 4, we introduce the FDM, as applied to direct and inverse problems, respectively. Numerical results are illustrated and discussed in Sections 5 and an extension of the study is presented in Section 6. Conclusions are provided in Section 7. 
\section{Mathematical Formulation}
The governing equation for a vibrating bounded structure $\Omega \subset \mathbb{R}^{n},\ n=1,2,3$, acted upon by a force $F(\underline{x},t)$ is given by the wave equation
\begin{eqnarray}
u_{tt}(\underline{x},t)=c^{2}\nabla^{2}u(\underline{x},t)+F(\underline{x},t),
\quad (\underline{x},t)\in\Omega\times(0,T), \label{eq1}
\end{eqnarray}
where $T>0$ is a given time, $u(\underline{x},t)$ represents the displacement and $c>0$ is the wave speed of propagation. For simplicity, we assume that $c$ is a constant, but we can also let $c$ be a function depending  on the space variable $\underline{x}$. For example, in $n=1$-dimension, where $\Omega$ represents the interval $(0,L), \ L>0$, occupied by a vibrating inhomogeneous string, its small transversal vibrations are governed by the wave equation
 \begin{eqnarray}
\omega(x)u_{tt}(x,t)=u_{xx}(x,t)+F(x,t),
\quad (x,t)\in(0,L)\times(0,T), \label{eqwx}
\end{eqnarray}
where $\omega(x)=c^{-2}(x)$ represents the mass density of the string, which is stretched by a unit force.

Equation \eqref{eq1} has to be solved subject to the initial conditions
\begin{equation}
u(\underline{x},0)=u_0(\underline{x}), \quad  \underline{x}\in \Omega,   \label{eq2}
\end{equation}
\begin{equation}
 u_{t}(\underline{x},0)=v_{0}(\underline{x}), \quad  \underline{x}\in \Omega, \label{eq3}
\end{equation}
where $u_0$ and $v_{0}$ represent the initial displacement and velocity, respectively. On the boundary of the structure $\partial{\Omega}$ we can prescribe Dirichlet, Neumann, Robin or mixed boundary conditions. 

Let us consider, for the sake of simplicity, Dirichlet boundary conditions being prescribed, namely, 
\begin{eqnarray}
u(\underline{x},t)=P(\underline{x},t), \quad  (\underline{x},t)\in\partial{\Omega}\times (0,T), \label{eq4}
\end{eqnarray}  
where $P$ is a prescribed boundary displacement.

If the force $F(\underline{x},t)$ is given, then equations \eqref{eq1}, \eqref{eq2}-\eqref{eq4} form a direct well-posed problem, see e.g. Morse and Feshbach (1953). However, if the force function $F(\underline{x},t)$ cannot be directly observed it hence becomes unknown and then clearly, the above set of equations is not sufficient to determine uniquely the pair solution $(u(\underline{x},t),F(\underline{x},t))$. Then, we consider the additional measurement of the flux tension of the structure on a (non-zero measure) portion $\Gamma\subset\partial{\Omega}$, namely,  
\begin{eqnarray}
\frac{\partial{u}}{\partial{\nu}}(\underline{x},t)=q(\underline{x},t), \quad (x,t)\in\Gamma\times(0,T), \label{eq6}
\end{eqnarray}
where $\underline{\nu}$ is the outward unit normal to $\partial{\Omega}$ and $q$ is a given function. Other additional information, such as the 'upper-base' final displacement measurement $u(\underline{x},T)$ for $\underline{x}\in\Omega$, will be investigated in a separate work.

Also, note that if instead of the Dirichlet boundary condition \eqref{eq4} we would have supplied a Neumann boundary condition then, the quantities $u$ and $\partial{u}/\partial{\nu}$ would have had to be reversed in \eqref{eq4} and \eqref{eq6}.  In order to ensure a unique solution we further assume that
\begin{eqnarray}
F(\underline{x},t)=f(\underline{x})h(\underline{x},t), \quad (\underline{x},t)\in \Omega\times(0,T), \label{eq8}
\end{eqnarray}
where $h(\underline{x},t)$ is a known function and $f(\underline{x})$ represents the unknown space-dependent forcing function to be determined. This restriction is necessary because otherwise, we can always add to $u(\underline{x},t)$ any function of the form $t^{2}U(\underline{x})$ with $U\in C^{2}(\overline{\Omega})$ arbitrary with compact support in $\Omega$, and still obtain another solution satisfying \eqref{eq1}, \eqref{eq2}-\eqref{eq6}.

Note that the unknown force $f(\underline{x})$ is an interior quantity and it depends on the space variable $\underline{x}\in \Omega \subset \mathbb{R}^{n}$, whilst the additional measurement \eqref{eq6} of the flux $q(\underline{x},t)$ is a boundary quantity and it depends on $(\underline{x},t)\in\Gamma\times(0,T)$. 

In the next subsection, we analyse more closely the uniqueness of solution of the inverse problem which requires finding the pair solution $(u(\underline{x},t),f(\underline{x}))$ satisfying equations \eqref{eq1}, \eqref{eq2}-\eqref{eq8}.
\subsection{Mathematical Analysis}
To start with, from \eqref{eq8}, and taking for simplicity $c=1$, equation \eqref{eq1} recasts as 
\begin{eqnarray}
u_{tt}(\underline{x},t)=
\nabla^{2}u(\underline{x},t)+f(\underline{x})h(\underline{x},t),
\quad (\underline{x},t)\in\Omega\times(0,T). \label{eq9}
\end{eqnarray}
We note that in the one-dimensional case, $n=1$, and for $c=h=1$ and other compatibility conditions satisfied by the data \eqref{eq2}-\eqref{eq6}, Cannon and Dunninger \cite{candunn70}, based on the method of Fourier series, established the uniqueness of a classical solution of the inverse problem. We also have the following more general uniqueness result, see Theorem 9 of \cite{engl94}.
\\
\\
{\bf Theorem 1.} {\it Assume that $\Omega \subset \mathbb{R}^{n}$ is a bounded star-shaped domain with sufficiently smooth boundary such that $T>diam(\Omega)$. Let $h\in H^{2}(0,T;L^{\infty}(\Omega))$ be such that $h(.,0)\in L^{\infty}(\Omega)$, $h_{t}(.,0)\in L^{\infty}(\Omega)$ and 
\begin{eqnarray}
H:=\frac{||h_{tt}||_{L^{2}(0,T;L^{\infty}(\Omega))}}{inf_{\underline{x}\in\Omega}|h(\underline{x},0)|} \quad \text{is sufficiently small}. \label{eqth1.1}
\end{eqnarray}
If $\Gamma=\partial\Omega$, then the inverse problem \eqref{eq2}-\eqref{eq6} and \eqref{eq9} has at most one solution $(u(\underline{x},t),f(\underline{x}))$ in the class of functions
\begin{eqnarray}
u\in L^{2}(0,T;H^{1}(\Omega)), \quad u_{t}\in L^{2}(0,T;L^{2}(\Omega)), \quad u_{tt}\in L^{2}(0,T;(H^{1}(\Omega))^{\prime}), \quad f\in L^{2}(\Omega), \label{eqth1.2}
\end{eqnarray}
where $(H^{1}(\Omega))^{\prime}$ denotes the dual of $H^{1}(\Omega)$.}\\
\\
For the notations and definitions of the function spaces involved, see \cite{lionsnewyork}.

The proof in \cite{engl94} relies on the estimate (5.25) of \cite{lions88}, namely,
\begin{eqnarray}
||h(.,0)f||_{L^{2}(\Omega)}\leq K_{1}||w_{1}||_{L^{2}(\partial{\Omega}\times(0,T))}, \label{eqth1.3}
\end{eqnarray}
for some positive constant $K_{1}$ which depends only on $\Omega$ and $T$, and $w_{1}$ is the solution of the problem
\begin{eqnarray}
w_{1tt}(\underline{x},t)=\nabla^{2}w_{1}(\underline{x},t), \quad (\underline{x},t)\in \Omega\times(0,T), \label{eqth1.4}
\end{eqnarray}
\begin{eqnarray}
w_{1}(\underline{x},0)=h(\underline{x},0)f(\underline{x}), \quad w_{1t}(\underline{x},0)=h_{t}(\underline{x},0)f(\underline{x}), \quad \underline{x}\in \Omega. \label{eqth1.6}
\end{eqnarray}
\begin{eqnarray}
\frac{\partial{w_{1}}}{\partial{\nu}}(\underline{x},t)=0, \quad (\underline{x},t)\in \partial{\Omega}\times(0,T), \label{eqth1.5}
\end{eqnarray}
Theorem 1 also requires that the quantity $H$ in equation \eqref{eqth1.1} is sufficiently small which can be guaranteed if $||h_{tt}||_{L^{2}(0,T;L^{\infty}(\Omega))}$ is small or, if $inf_{x\in \Omega}|h(x,0)|$ is large. For example, if
\begin{eqnarray}
h(\underline{x},t)=t \ h_{1}(\underline{x})+h_{2}(\underline{x}), \quad (\underline{x},t)\in \Omega\times(0,T), \label{eqth1.7}
\end{eqnarray}
with $h_{1}\in L^{\infty}(\Omega)$ and $h_{2}\in L^{\infty}(\Omega)$ given functions, then $h_{tt}=0$ and therefore condition \eqref{eqth1.1} is satisfied if  $inf_{x\in \Omega}|h_{2}(\underline{x})|>0$. In this case, the uniqueness proof follows immediately by remarking that $w_{1}=u_{tt}$, where $u$ satisfies the problem
\begin{eqnarray}
u_{tt}=\nabla^{2}u(\underline{x},t)+(t \ h_{1}(\underline{x})+h_{2}(\underline{x}))f(\underline{x}), \quad (\underline{x},t)\in \Omega\times(0,T), \label{eqth1.8}
\end{eqnarray}
\begin{eqnarray}
u(\underline{x},0)=u_{t}(\underline{x},0)=0, \quad \underline{x}\in \Omega, \label{eqth1.9}
\end{eqnarray}
\begin{eqnarray}
\frac{\partial{u}}{\partial{\nu}}(\underline{x},0)=0, \quad (\underline{x},t)\in \partial{\Omega}\times(0,T). \label{eqth1.10}
\end{eqnarray}
In the above, $(u,f)$ represents the difference between two solutions $(u_{1},f_{1})$ and $(u_{2},f_{2})$ of the inverse problem \eqref{eq2}-\eqref{eq6} and \eqref{eq9}. Then from \eqref{eq4} it follows that
\begin{eqnarray}
u(\underline{x},t)=0, \quad (\underline{x},t)\in \partial{\Omega}\times(0,T). \label{eqth1.11}
\end{eqnarray}
Since $w_{1}=u_{tt}$, from \eqref{eqth1.11} it results that 
\begin{eqnarray}
w_{1}(\underline{x},t)=0, \quad (\underline{x},t)\in \partial{\Omega}\times(0,T). \label{eqth1.12}
\end{eqnarray}
Then, conditions on $\Omega$ and $T>diam(\Omega)$, and equations \eqref{eqth1.4}, \eqref{eqth1.5} and \eqref{eqth1.12}, implies that the uniqueness property, see Remark 1.7 of \cite{lions88}, is applicable and consequently, $w_{1}\equiv0$.  
Then, \eqref{eqth1.3} and  $inf_{x\in \Omega}|h_{2}(\underline{x})|>0$ immediately yields $f\equiv0$. Afterwords, the problem \eqref{eqth1.8}-\eqref{eqth1.10} with $f=0$ yields $u\equiv0$.
\\
\\
We also have the following uniqueness theorem due to Theorem 3.8 of \cite{vmk92}.
\\
\\
{\bf Theorem 2.} {\it Assume that $\Omega\subset\mathbb{R}^{n}$ is a bounded domain with piecewise smooth boundary. Let $h\in C^{3}(\overline{\Omega}\times[0,T])$ be such that
\begin{eqnarray}
h(\underline{x},0)\neq0 \quad \text{for} \quad \underline{x}\in\overline{\Omega}. \label{eqth2.1}
\end{eqnarray}
If $\Gamma=\partial{\Omega}$, then the inverse problem \eqref{eq2}-\eqref{eq6} and \eqref{eq9} has at most one solution $(u(\underline{x},t),f(\underline{x}))\in C^{3}(\overline{\Omega}\times[0,T])\times C(\overline{\Omega})$.}\\
\\
One can remark that the previously stated uniqueness Theorems 1 and 2 require that the Neumann observation \eqref{eq6} is over the complete boundary $\Gamma=\partial{\Omega}$. In the incomplete case that $\Gamma\subset\partial{\Omega}$ is only a part of $\partial{\Omega}$ then, the uniqueness Theorem 1 holds under the assumption that $h$ is independent of $\underline{x}$, \cite{my95}, as follows.\\
\\
{\bf Theorem 3.} {\it Assume that $\Omega\subset\mathbb{R}^{n}$ is a bounded star-shaped domain with smooth boundary such that $T>diam(\Omega)$. Let $h\in C^{1}[0,T]$ be independent of $\underline{x}$ such that equation \eqref{eq9} becomes
\begin{eqnarray}
u_{tt}(\underline{x},t)=\nabla^{2}u(\underline{x},t)+f(\underline{x})h(t), \quad (\underline{x},t)\in \Omega\times(0,T), \label{eqth3.1}
\end{eqnarray}
and assume further that $h(0)\neq0$. Then the inverse problem \eqref{eq2}-\eqref{eq6} and \eqref{eqth3.1} has at most one solution in the class of functions}
\begin{eqnarray}
u\in C^{1}([0,T];H^{1}(\Omega))\cap C^{2}([0,T];L^{2}(\Omega)), \quad f\in L^{2}(\Omega). \label{eqth3.2}
\end{eqnarray}

In Section 4, we shall consider the numerical determination of space-dependent forcing functions. But before we do that, in the next section we explain the finite-difference method (FDM) adopted for the numerical discretisation of the direct problem.
\section{Numerical Solution of the Direct Problem}
In this section, we consider the direct initial Dirichlet boundary value problem \eqref{eq1}, \eqref{eq2}-\eqref{eq4} for simplicity, in one-dimension, i.e. $n=1$ and $\Omega=(0,L)$ with $L>0$, when the force $F(x,t)$
is known and the displacement $u(x,t)$ is to be determined, namely,  
\begin{eqnarray}
u_{tt}(x,t)=c^{2}u_{xx}(x,t)+F(x,t),
\quad (x,t)\in(0,L)\times(0,T], \label{eq11}
\end{eqnarray}
\begin{eqnarray}
u(x,0)=u_0(x), \quad  u_{t}(x,0)=v_{0}(x), \quad  x\in[0,L], \   \label{eq12}
\end{eqnarray}
\begin{eqnarray}
u(0,t)=P(0,t)=:P_{0}(t), \quad t\in(0,T],  \label{eq13}
\end{eqnarray}  
\begin{eqnarray}
u(L,t)=P(L,t)=:P_{L}(t), \quad  t\in(0,T].    \label{eq14}
\end{eqnarray}
The compatibility conditions between \eqref{eq12}-\eqref{eq14} yield
\begin{eqnarray}
P_{0}(0)=u_{0}(0), \quad P_{L}(0)=u_{0}(L).  \label{eq15}
\end{eqnarray}

The discrete form of our problem is as follows. We divide the domain $(0,L)\times(0,T)$ into $M$ and $N$ subintervals of equal space
length $\delta x$ and $\delta t$, where  $\delta x=L/M$ and $\delta t=T/N$. We denote by $u_{i,j}:=u(x_{i},t_{j})$, where $x_{i}=i\delta x$, $t_{j}=j\delta t$, and $F_{i,j}:=F(x_{i},t_{j})$
for $i=\overline{0,M}$, $j=\overline{0,N}$. Then, a central-difference approximation to equations \eqref{eq11}-\eqref{eq14} at the mesh points $(x_{i},t_{j})=(i\delta x,j\delta t)$ 
of the rectangular mesh covering the solution domain $(0,L)\times(0,T)$ is, \cite{gds85}, 
\begin{eqnarray}
u_{i,j+1}=r^{2}u_{i+1,j}+2(1-r^{2})u_{i,j}+r^{2}u_{i-1,j}-u_{i,j-1}+(\delta t)^{2}F_{i,j}, \label{eq16}\\
\quad \quad  \quad i=\overline{1,(M-1)}, \quad j=\overline{1,(N-1)},\notag 
\end{eqnarray}
\begin{eqnarray}
u_{i,0}=u_{0}(x_{i}), \quad i=\overline{0,M}, \quad \frac{u_{i,1}-u_{i,-1}}{2(\delta t)}=v_{0}(x_{i}), \quad i=\overline{1,(M-1)}, \label{eq16.1}
\end{eqnarray}
\begin{eqnarray}
u_{0,j}=P_{0}(t_{j}), \quad u_{M,j}=P_{L}(t_{j}), \quad j=\overline{0,N},\label{eq16.2}
\end{eqnarray}
where $r=c(\delta t)/\delta x$. Equation \eqref{eq16} represents an explicit FDM which is stable if $r\leq1$,  giving approximate values for the solution at mesh points along
$t=2\delta t,3\delta t,...,$ as soon as the solution at the mesh points along $t=\delta t$ has been determined. Putting 
$j=0$ in equation \eqref{eq16} and using \eqref{eq16.1}, we obtain
\begin{eqnarray}
u_{i,1}=\frac{1}{2}r^{2}u_{0}(x_{i+1})+(1-r^{2})u_{0}(x_{i})+\frac{1}{2}r^{2}u_{0}(x_{i-1})+(\delta t)v_{0}(x_{i})+\frac{1}{2}(\delta t)^{2}F_{i,0}, \notag \\
\quad \quad \quad \quad \quad i=\overline{1,(M-1)}. \label{eq16.3}
\end{eqnarray}  
The normal derivatives $\frac{\partial{u}}{\partial{\nu}}(0,t)$ and 
$\frac{\partial{u}}{\partial{\nu}}(L,t)$ are calculated using the finite-difference approximations 
\begin{eqnarray}
-\frac{\partial{u}}{\partial{x}}(0,t_{j})=-\frac{4u_{1,j}-u_{2,j}-3u_{0,j}}{2(\delta x)}, \quad \quad 
\frac{\partial{u}}{\partial{x}}(L,t_{j})=\frac{3u_{M,j}-4u_{M-1,j}+u_{M-2,j}}{2(\delta x)},  \quad \notag \\ j=\overline{1,N}. \label{eq19}
\end{eqnarray}
\section{Numerical Solution of the Inverse Problem}
We now consider the inverse initial boundary value problem \eqref{eq2}-\eqref{eq6}  and \eqref{eq9} in one-dimension, i.e. $n=1$ and $\Omega=(0,L)$, when both the force $f(x)$ and the displacement $u(x,t)$ are to be determined, from the governing equation (take $c=1$ for simplicity)
\begin{eqnarray}
u_{tt}(x,t)=u_{xx}(x,t)+f(x)h(x,t), \quad (x,t)\in(0,L)\times(0,T], \label{eqfsplite}
\end{eqnarray}
 subject to the initial and boundary conditions \eqref{eq12}-\eqref{eq14} and the overspecified flux tension condition \eqref{eq6} at one end of the string, say at $x=0$, namely
 \begin{eqnarray}
 -\frac{\partial{u}}{\partial{x}}(0,t)=q(0,t)=:q_{0}(t), \quad t\in(0,T]. \label{equofxat0t}
 \end{eqnarray}
 
In the case that $h$ is independent of $x$, according to Theorem 3, the inverse source problem \eqref{eq12}-\eqref{eq14}, \eqref{eqfsplite} and \eqref{equofxat0t} has at most one solution provided that $h\in C^{1}[0,T]$, $h(0)\neq0$ and $T>L$.

In discretised finite-difference form equations \eqref{eq12}-\eqref{eq14} and \eqref{eqfsplite}  recast as equations \eqref{eq16.1}, \eqref{eq16.2},
\begin{eqnarray}
u_{i,j+1}-(\delta t)^{2}f_{i}h_{i,j}=r^{2}u_{i+1,j}+2(1-r^{2})u_{i,j}+r^{2}u_{i-1,j}-u_{i,j-1}, \label{eqfdmsplitef}\\
\quad \quad  \quad i=\overline{1,(M-1)}, \quad j=\overline{1,(N-1)},\notag 
\end{eqnarray}
and
\begin{eqnarray}
u_{i,1}-\frac{1}{2}(\delta t)^{2}f_{i}h_{i,0}=\frac{1}{2}r^{2}u_{0}(x_{i+1})+(1-r^{2})u_{0}(x_{i})+\frac{1}{2}r^{2}u_{0}(x_{i-1})+(\delta t)v_{0}(x_{i}), \label{eqatjiszero}\\
\quad \quad  \quad i=\overline{1,(M-1)},\notag 
\end{eqnarray}
where $f_{i}:=f(x_{i})$ and $h_{i,j}:=h(x_{i},t_{j})$.

Discretizing \eqref{equofxat0t} using \eqref{eq19} we also have
\begin{eqnarray}
q_{0}(t_{j})=-\frac{\partial{u}}{\partial{x}}(0,t_{j})=-\frac{4u_{1,j}-u_{2,j}-3u_{0,j}}{2(\delta x)},\quad  j=\overline{1,N}. \label{equx0}
\end{eqnarray}

In practice, the additional observation \eqref{equx0} comes from measurement which is inherently contaminated with errors. We therefore model this by replacing the exact data $q_{0}(t)$ by the noisy data 
\begin{eqnarray}
q_{0}^{\epsilon}(t_{j})=q_{0}(t_{j})+\epsilon_{j}, \ \ \ j=\overline{1,N},\label{eq4.5}
\end{eqnarray}
where $(\epsilon_{j})_{j=\overline{1,N}}$ are $N$ random noisy variables generated (using the MATLAB  routine 'normrd') from a Gaussian normal distribution with mean zero and standard deviation $\sigma=p\times max_{t\in[0,T]}\left|q_{0}(t)\right|$, where $p$ represents the percentage of noise. 

Assembling \eqref{eqfdmsplitef}-\eqref{equx0} and using \eqref{eq16.1} and \eqref{eq16.2}, the discretised inverse problem reduces to solving a global linear system of $(M-1)\times N+N$ equations with $(M-1)\times N+(M-1)$ unknowns. Since this system is linear we can eliminate the unknowns $u_{i,j}$ for $i=\overline{1,(M-1)}$, $j=\overline{1,N}$, to reduce the problem to solving an ill-conditioned system of $N$ equations with $(M-1)$ unknowns of the generic form  
\begin{eqnarray}
A\underline{f}=\underline{b}^{\epsilon}, \label{eqA}
\end{eqnarray}
where the right-hand side vector $\underline{b}^{\epsilon}$ incorporates the noisy measurement \eqref{eq4.5}. For a unique solution we require $N\geq M-1$. The method of least squares can be used to find an approximate solution to overdetermined systems. For the system of equations \eqref{eqA}, the least squares solution is given by
$\underline{f} =(A^{tr}A)^{-1}A^{tr}\underline{b}^{\epsilon}$, where the superscript $^{tr}$ denotes the transpose. 

\begin{figure}[H]
    \begin{subfigure}[b]{0.3\textwidth}
     \caption{}
        \centering
        \includegraphics[width=11cm]{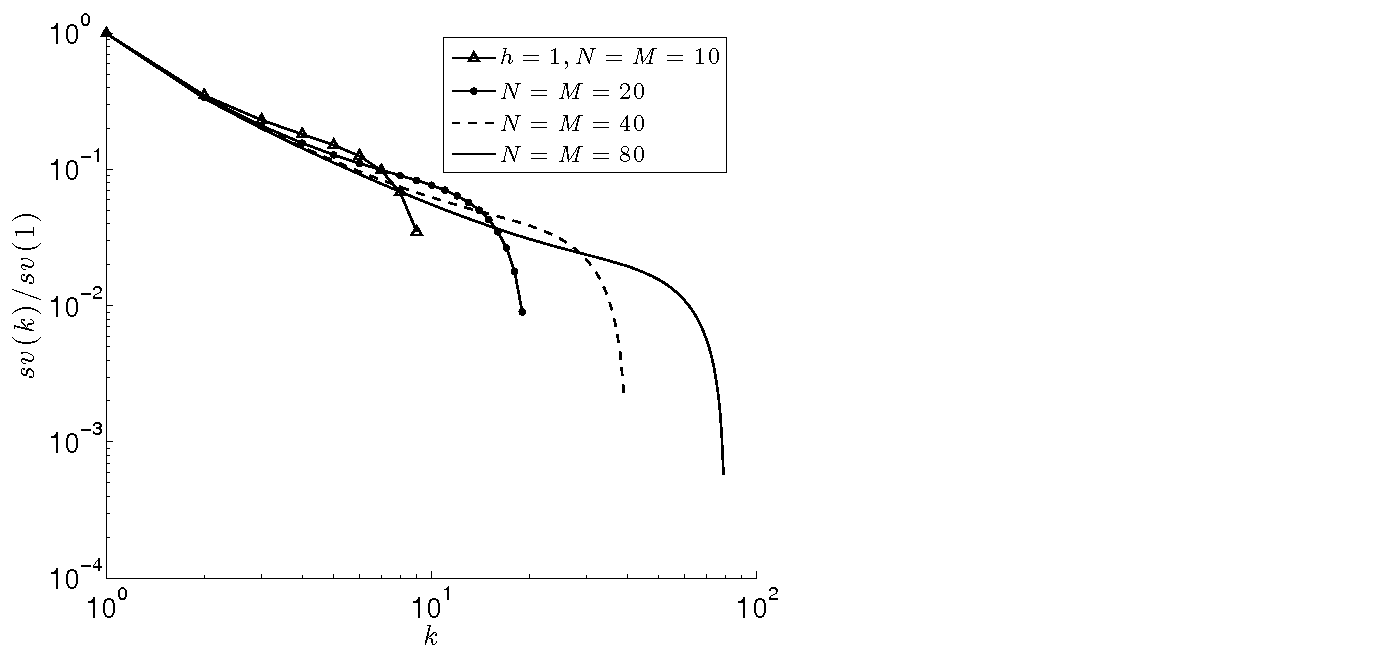}
    \end{subfigure}
    \hspace{3.3cm}
    \begin{subfigure}[b]{0.3\textwidth}
     \caption{}
        \centering
        \includegraphics[width=11cm]{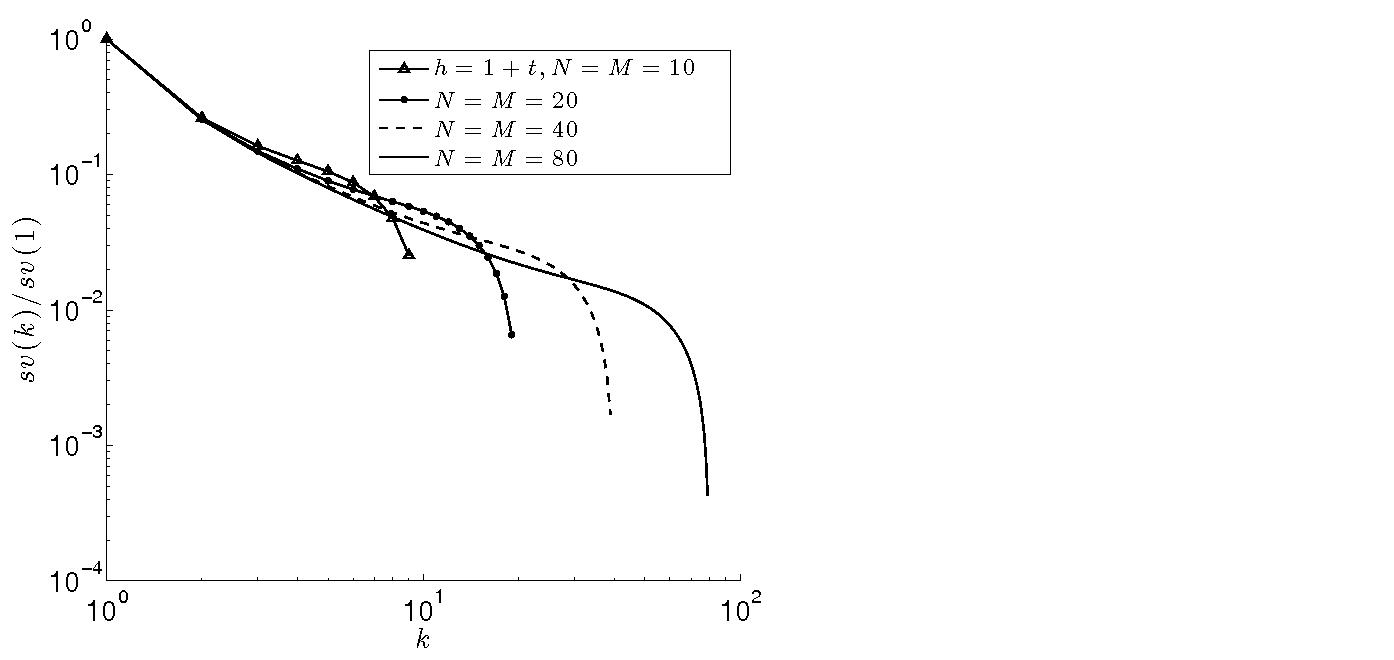}
            \end{subfigure}\\
    \begin{subfigure}[b]{0.3\textwidth}
     \caption{}
        \includegraphics[width=11cm]{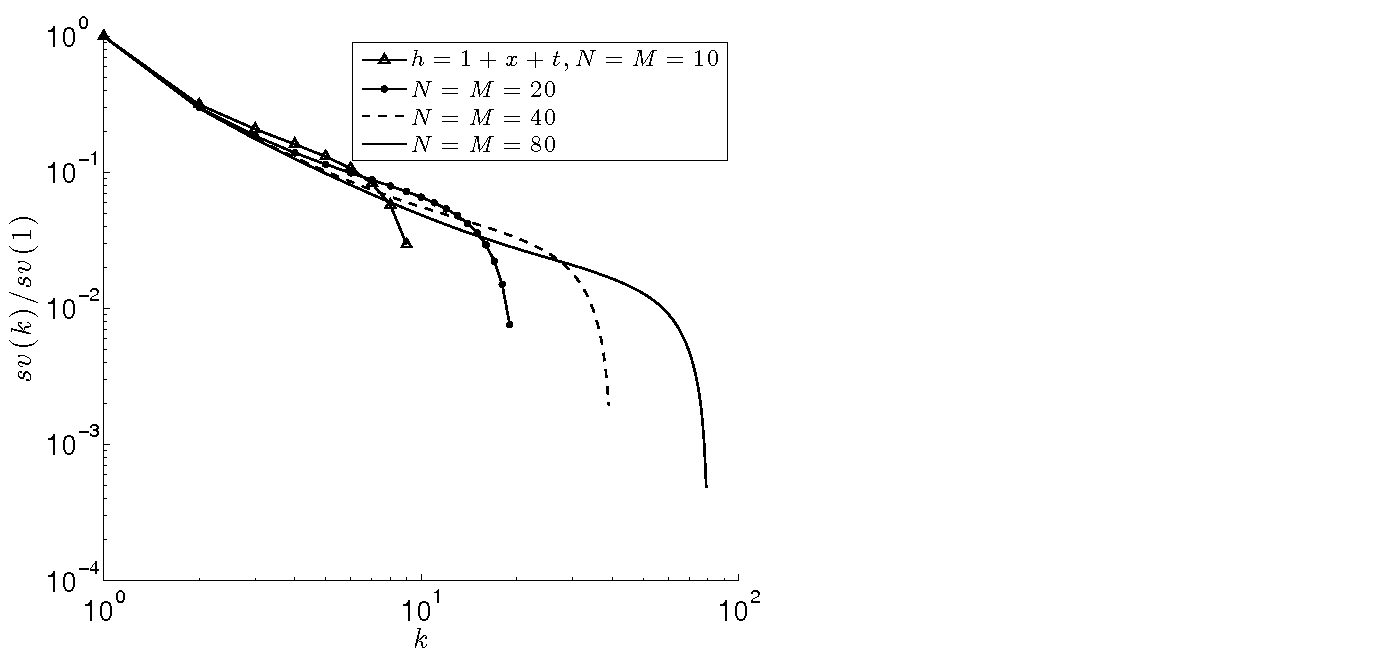}
    \end{subfigure}
    \hspace{3.3cm}
    \begin{subfigure}[b]{0.3\textwidth}
     \caption{}
        \includegraphics[width=11cm]{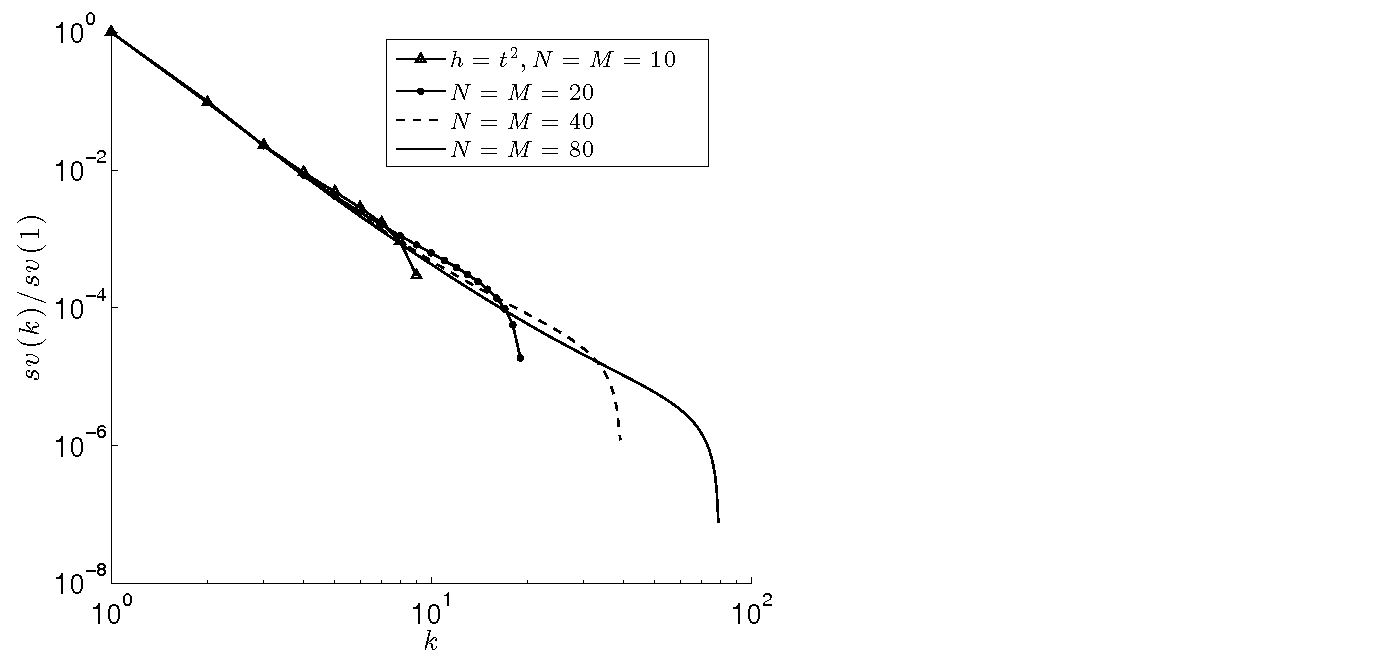}
    \end{subfigure}
    \caption{Normalised singular values $sv(k)/sv(1)$ for $k=\overline{1,(M-1)}$, for (a) Example 1, (b) Example 2, (c) Example 3, and (d) Example 4.}
    \label{fig:normizedsvdEx1Ex2Ex3Ex4problem3}
\end{figure}

For the Examples 1-4 that will be considered in the next section, the condition numbers of the matrix $A$ in \eqref{eqA} (calculated using the command cond($A$) in MATLAB) given in Table 1 are between O($10^{4}$) to O($10^{8}$) for $M=N=80$. These large condition numbers indicate that the system of equations \eqref{eqA} is ill-conditioned. The ill-conditioning nature of the matrix $A$ can also be revealed by plotting its normalised singular values $sv(k)/sv(1)$ for $k=\overline{1,(M-1)}$, in Figure \ref{fig:normizedsvdEx1Ex2Ex3Ex4problem3}. These singular values have been calculated in MATLAB using the command svd($A$).
\begin{table}[H]
\caption{Condition number of matrix $A$ for Examples 1-4.}
\centering
\begin{tabular}{|c|c|c|c|c|c}
\hline
$$ & Example 1 & Example 2 & Example 3 & Example 4 \\
$N=M$ & $h(x,t)=1$ &  $h(x,t)=1+t$  & $h(x,t)=1+x+t$ & $h(x,t)=t^{2}$ \\
\hline
$10$	 &$28.55$   &$39.53$ & $33.73$ &$3394.55$ \\
\hline
$20$	 &$110.98$  &$152.38$  & $131.29$ &$53232.36$ \\
\hline
$40$	 &$437.93$  &$596.91$  & $518.51$ &$826827.12$ \\
\hline
$80$	 &$1740.25$  &$2361.22$  & $2061.53$ &$12956244.4$ \\
\hline
\end{tabular}
\end{table}
\section{Numerical Results and Discussion}
In all examples in this section we take, for simplicity, $c=L=T=1$. Although the geometrical condition $1=T>diam(\Omega)=L=1$ is slightly violated, it is expected that the uniqueness Theorems 1 and 3 still hold, especially in $n=1$-dimension and when the inverse problems are numerically discretised.
\subsection{Example 1 ($h(x,t)=1$)}
This is an example in which we take $h(x,t)=1$ a constant function and consider first the direct problem \eqref{eq11}-\eqref{eq14}  with the input data
\begin{eqnarray}
u(x,0)=u_{0}(x)=\sin(\pi x), \quad
u_{t}(x,0)=v_{0}(x)=1, \quad x\in[0,1],   \label{eq20}
\end{eqnarray}
\begin{eqnarray}
u(0,t)=P_{0}(t)=t+\frac{t^{2}}{2},\quad u(1,t)=P_{L}(t)=t+\frac{t^{2}}{2}, \quad t\in(0,1], \label{eq21}
\end{eqnarray}
\begin{eqnarray}
F(x,t)=f(x)=1+\pi^{2}\sin(\pi x), \ \ \  x\in(0,1). \label{eq23}
\end{eqnarray}
The exact solution is given by
\begin{eqnarray}
u(x,t)&=&\sin(\pi x)+t+\frac{t^{2}}{2},\ \ \ (x,t)\in[0,1]\times[0,1].  \label{eq22}
\end{eqnarray}
  
\begin{figure}[H]
\begin{center}
\includegraphics[width=15cm]{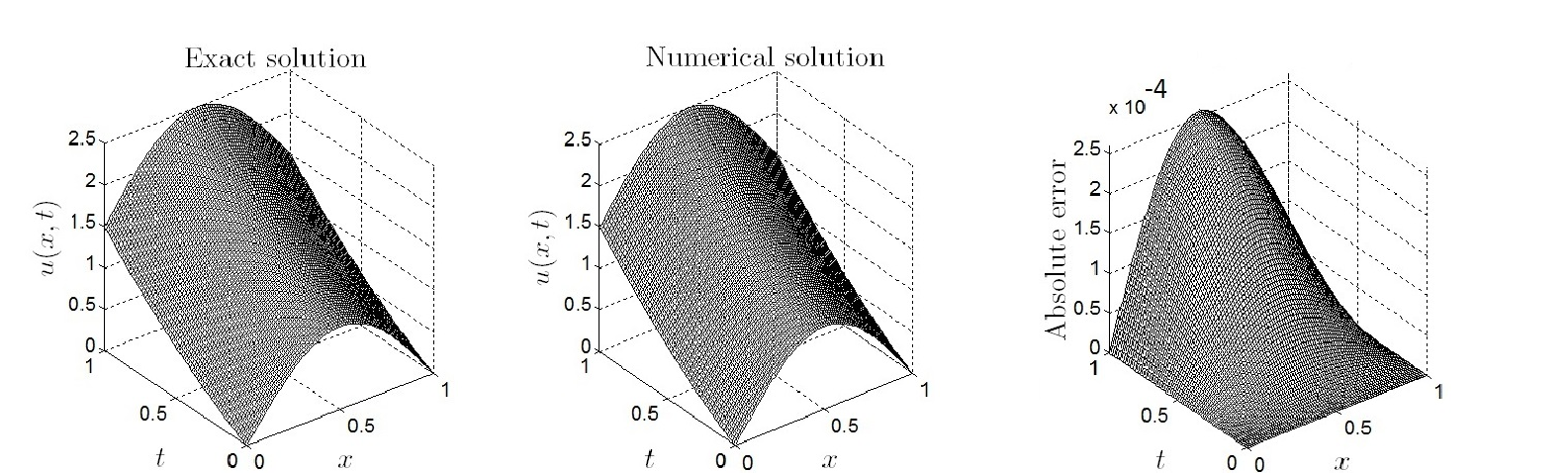}
\caption{Exact and numerical solutions for the displacement $u(x,t)$ and the absolute error between them for the direct problem obtained with $N=M=80$, for Example 1.}
\label{figexactnumericalerrorofuEx1problem3}
\end{center}
\end{figure}

The numerical and exact solutions for $u(x,t)$ at interior points are shown in Figure \ref{figexactnumericalerrorofuEx1problem3} and one can observe that an excellent agreement is obtained. Table 2 also gives the exact and numerical solutions for the flux tension \eqref{equofxat0t}. From this table it can be seen that the numerical results are convergent, as the mesh size decreases, and they are in very good agreement with the exact solution \eqref{eqq0}. Although not illustrated, it is reported that the same excellent agreement has been obtained between the exact and numerical solutions for the flux tension at $x=1$ and therefore, they are not presented. 
\begin{table}[H]
\caption{Exact and numerical solutions for the flux tension at $x=0$, for the direct problem of Example 1.}
\centering
\begin{tabular}{|c|c|c|c|c|c|c|}
\hline
$t$  &$0.1$   &$0.2$  &$...$   &$0.8$ &$0.9$ &$1$ \\  
\hline
$N=M=10$	&$-3.2427$  &$-3.2465$  &$...$	&$-3.2899$	 &$-3.2937$  &$-3.295$\\
\hline
$N=M=20$	&$-3.1675$  &$-3.1685$  &$...$	&$-3.1790$	 &$-3.1799$  &$-3.1802$\\
\hline
$N=M=40$	&$-3.1481$  &$-3.1483$  &$...$	&$-3.1510$	 &$-3.1512$  &$-3.1513$\\
\hline
$N=M=80$	&$-3.1432$  &$-3.1433$  &$...$	&$-3.1439$	 &$-3.1440$   &$-3.1440$\\
\hline
$exact$	&$-3.1416$  &$-3.1416$  &$...$	&$-3.1416$	 &$-3.1416$   &$-3.1416$\\
\hline
\end{tabular}
\end{table}
The inverse problem given by equations \eqref{eqfsplite} with $h(x,t)=1$, \eqref{eq20}, \eqref{eq21} and 
\begin{eqnarray}
-\frac{\partial{u}}{\partial{x}}(0,t)=q_{0}(t)=-\pi, \quad t\in[0,1], \label{eqq0}
\end{eqnarray}
is considered next. Since $h(0)=1\neq0$, Theorem 3 ensures the uniqueness of the solution in the class of functions \eqref{eqth3.2}, which in $n=1$-dimension rewrites as
\begin{eqnarray}
u\in C^{1}([0,T];H^{1}(0,L))\cap C^{2}([0,T];L^{2}(0,L)), \quad f\in L^{2}(0,L). \label{eqth3.2.1}
\end{eqnarray}
 In fact, the exact solution $(f(x),u(x,t))$ of this inverse problem is given by equations \eqref{eq23} and \eqref{eq22}, respectively. Numerically, we employ the FDM for discretising the inverse problem, as described in Section 4.
\subsubsection{Exact Data}
\begin{figure}[H]
    \begin{subfigure}[b]{0.3\textwidth}
     \caption{}
        \centering
        \includegraphics[width=11cm]{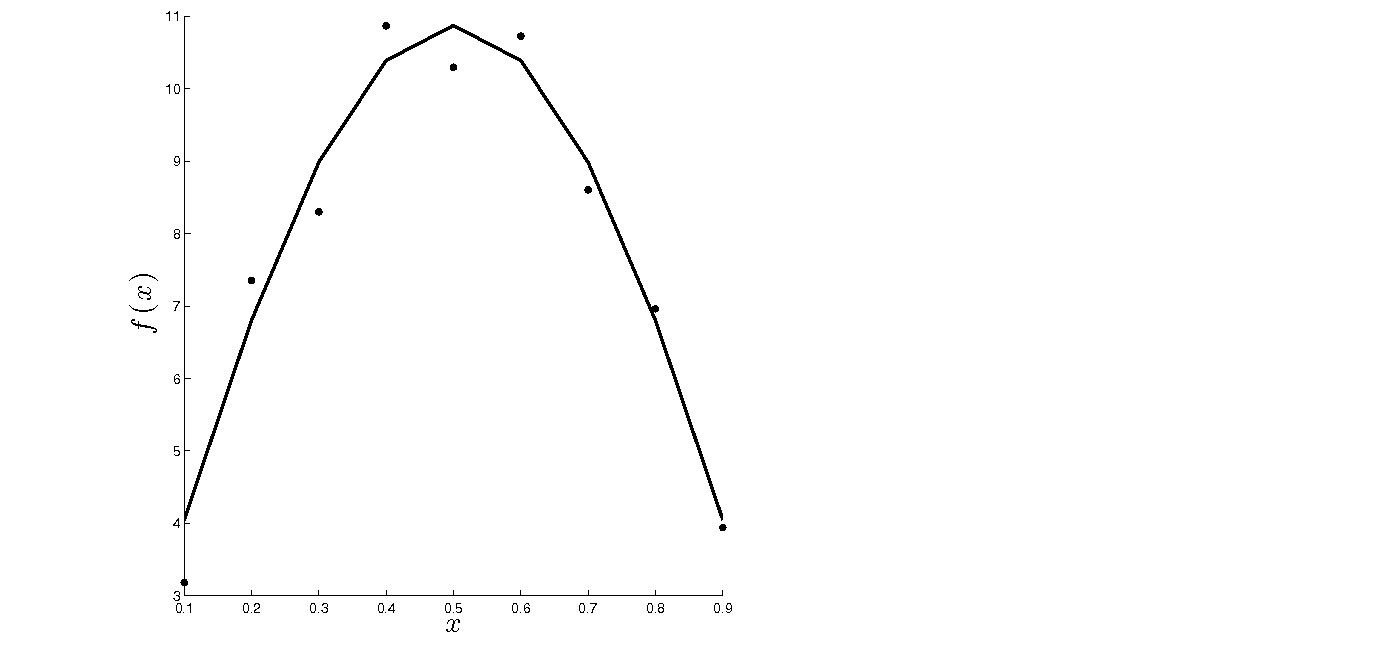}
    \end{subfigure}
    \hspace{2.9cm}
    \begin{subfigure}[b]{0.3\textwidth}
     \caption{}
        \centering
        \includegraphics[width=11cm]{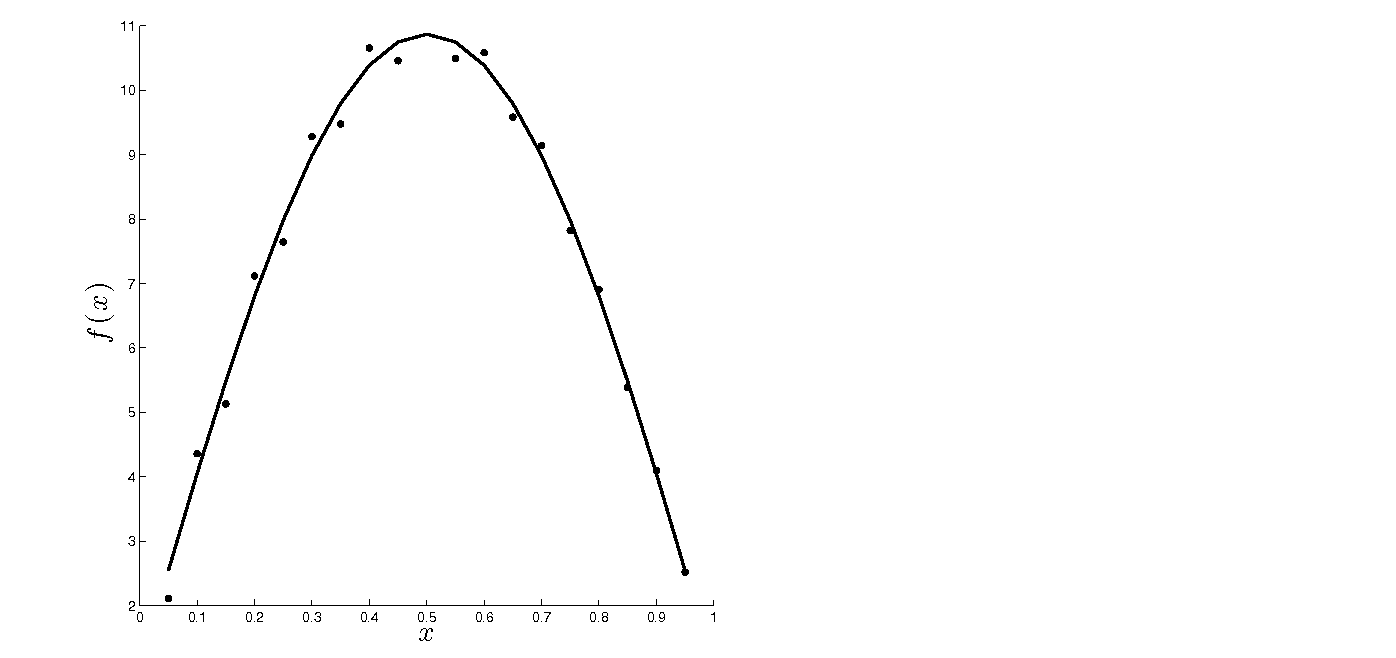}
    \end{subfigure}\\
\\
    \begin{subfigure}[b]{0.3\textwidth}
     \caption{}
        \centering
        \includegraphics[width=11cm]{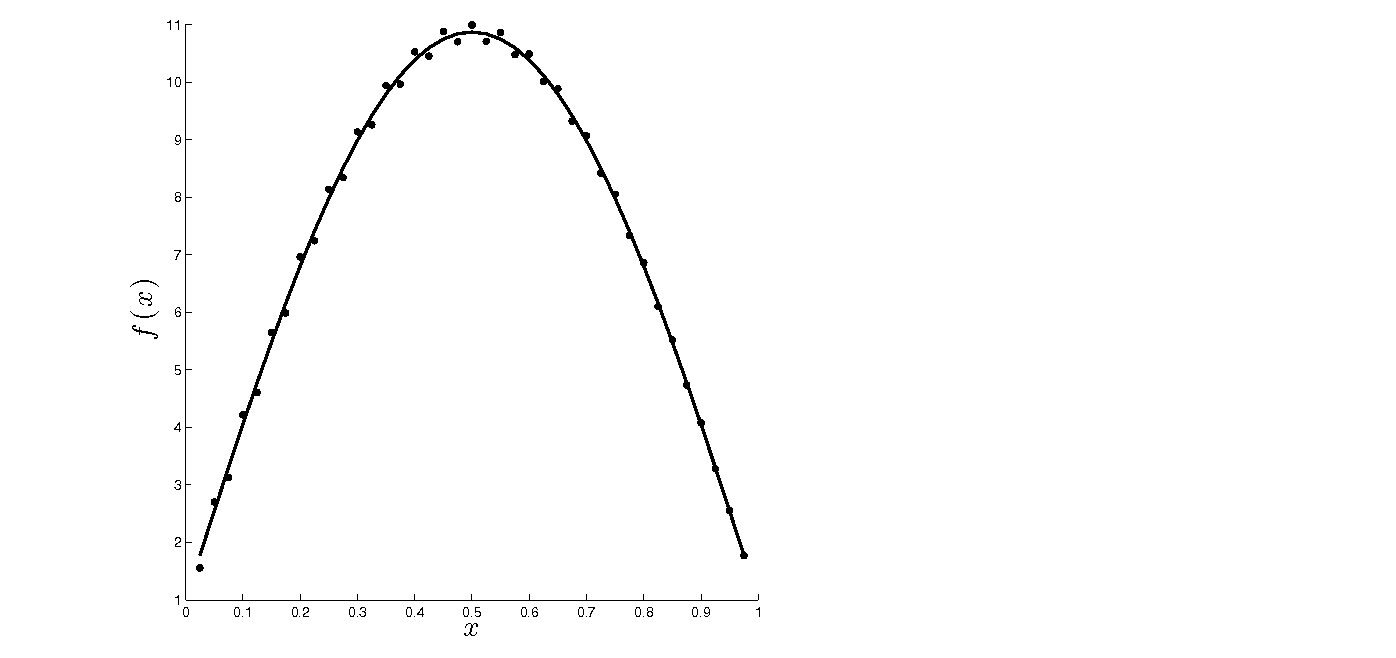}
    \end{subfigure}
\hspace{2.9cm}
    \begin{subfigure}[b]{0.3\textwidth}
     \caption{}
        \centering
        \includegraphics[width=11cm]{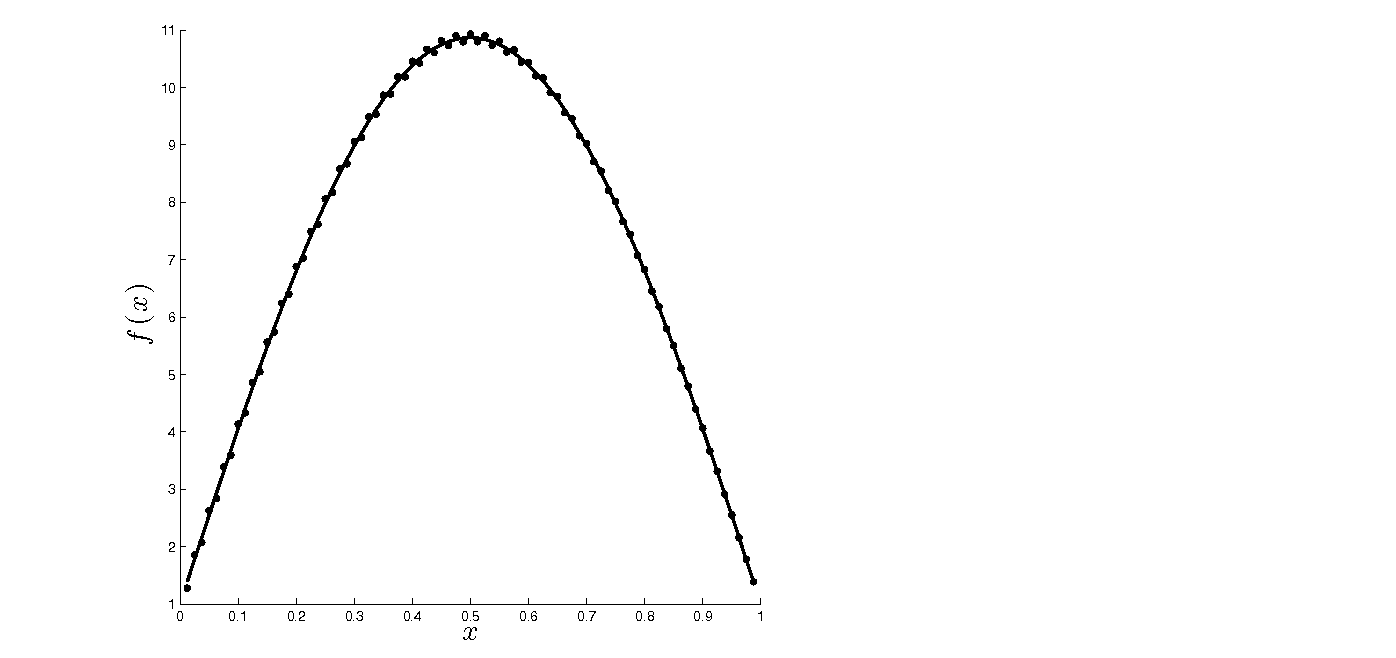}
    \end{subfigure}
    \caption{The exact (---) solution \eqref{eq23} for the force $f(x)$ in comparison with the numerical solution (\big{.\cdot\cdot\cdot}) for various $N=M=\text{(a)} \ 10, \text{(b)} \ 20, \text{(c)} \ 40,\text{and (d)} \ 80$,  and no regularization, for exact data, for the inverse problem of Example 1.}
    \label{fig:inverseproblemfexactnumericalEx1problem3}
\end{figure}

\begin{figure}[H]
    \begin{subfigure}[b]{0.3\textwidth}
    \caption{}
        \centering
        \includegraphics[width=13cm]{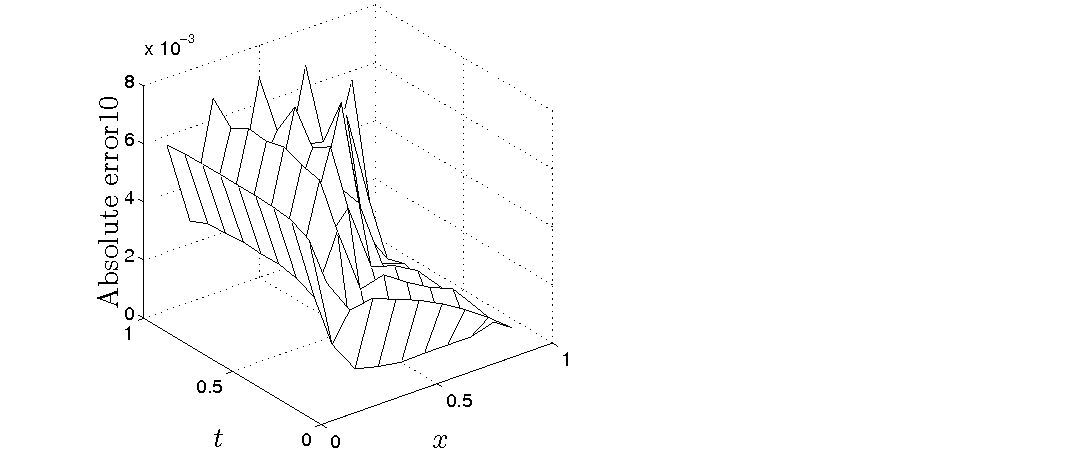}
    \end{subfigure}
    \hspace{2.9cm}
    \begin{subfigure}[b]{0.3\textwidth}
       \caption{}
        \centering
        \includegraphics[width=13cm]{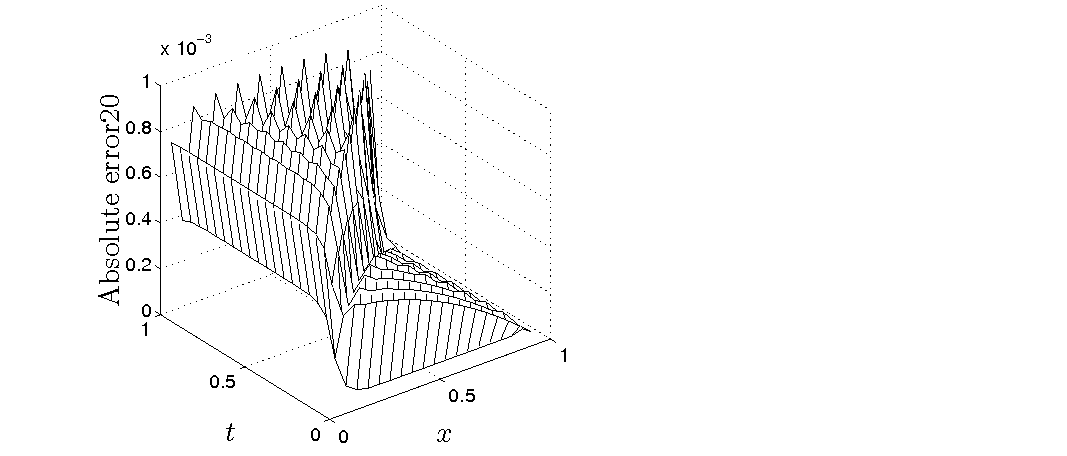}
    \end{subfigure}\\
\\
    \begin{subfigure}[b]{0.3\textwidth}
       \caption{}
        \centering
        \includegraphics[width=13cm]{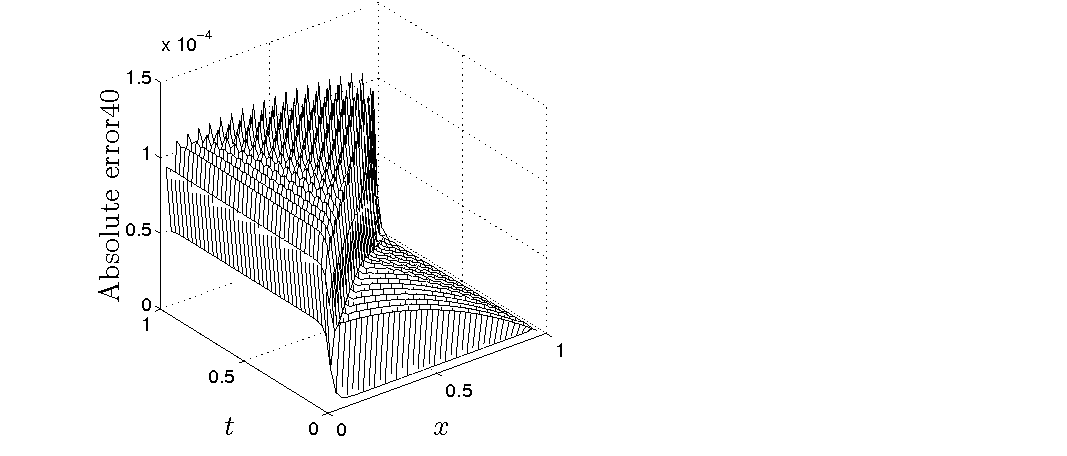}
    \end{subfigure}
\hspace{2.9cm}
    \begin{subfigure}[b]{0.3\textwidth}
      \caption{}
        \centering
        \includegraphics[width=13cm]{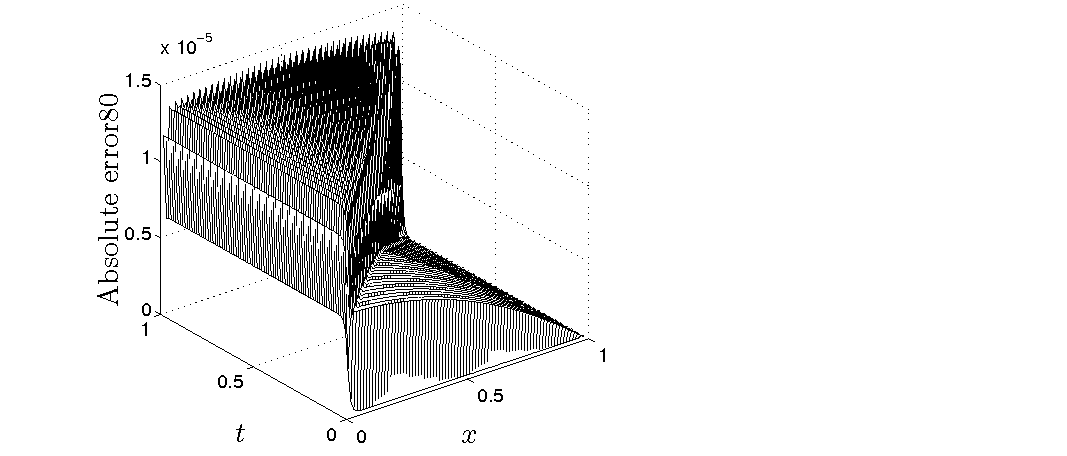}
    \end{subfigure}
    \caption{The absolute errors between the exact and numerical displacement $u(x,t)$ obtained with $N=M\in\lbrace10,20,40,80\rbrace$  and no regularization, for exact data, for the inverse problem of Example 1.}
    \label{fig:inverseproblemabsoulterroeofuEx1problem3}
\end{figure}

We first consider the case of exact data, i.e. $p=0$ and hence $\underline{\epsilon}=\underline{0}$ in \eqref{eq4.5}. The numerical results corresponding to $f(x)$ and $u(x,t)$ are plotted in Figures \ref{fig:inverseproblemfexactnumericalEx1problem3} and \ref{fig:inverseproblemabsoulterroeofuEx1problem3}, respectively. From these figures it can be seen that convergent and accurate numerical solutions are obtained.
\subsubsection{Noisy Data}
\begin{figure}[H]
        \begin{center}
        \includegraphics[width=15cm]{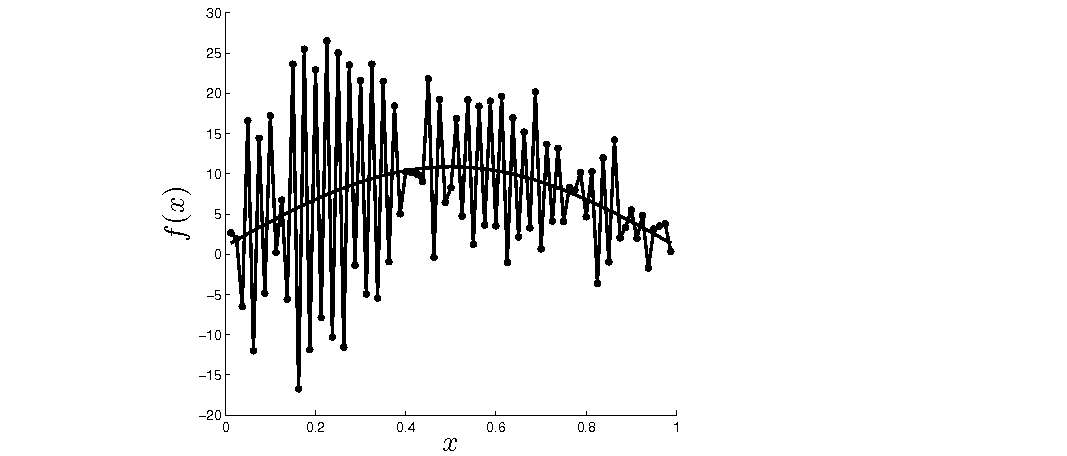}
\caption{The exact solution \eqref{eq23} for the force $f(x)$ in comparison with the numerical solution (\big{.\cdot\cdot\cdot}) for $N=M=80$, with no regularization, for $p=1\%$ noisy data, for the inverse problem of Example 1.}
\label{fig:inverseproblemaddnoiseatlambda0Ex1problem3}
\end{center}
\end{figure}

\begin{figure}[H]
\begin{center}
\includegraphics[width=17cm]{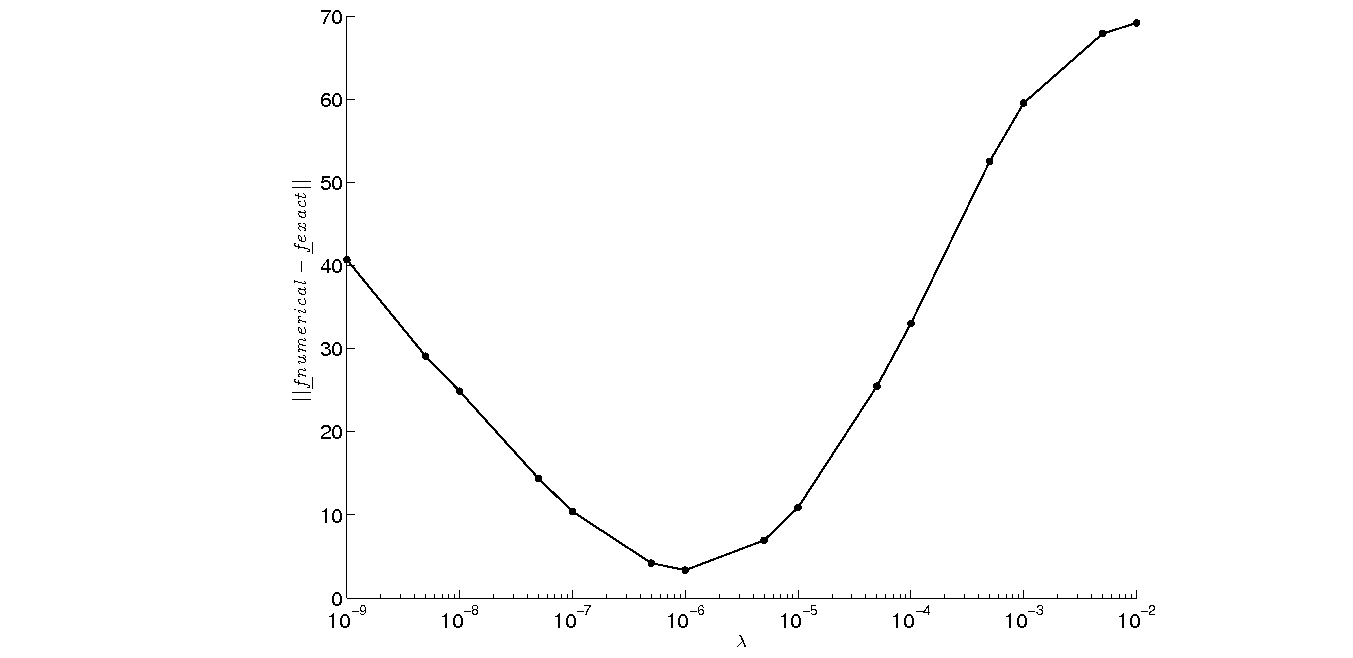}
\end{center}
\caption{The accuracy error $||\underline{f}numerical-\underline{f}exact||$, as a function of $\lambda$, for $N=M=80$ and $p=1\%$ noise, for the inverse problem of Example 1.}
\label{fig:normoffatexactandnumericalEx1problem3}
\end{figure}

\begin{figure}[H]
       \begin{center}
      \includegraphics[width=17cm]{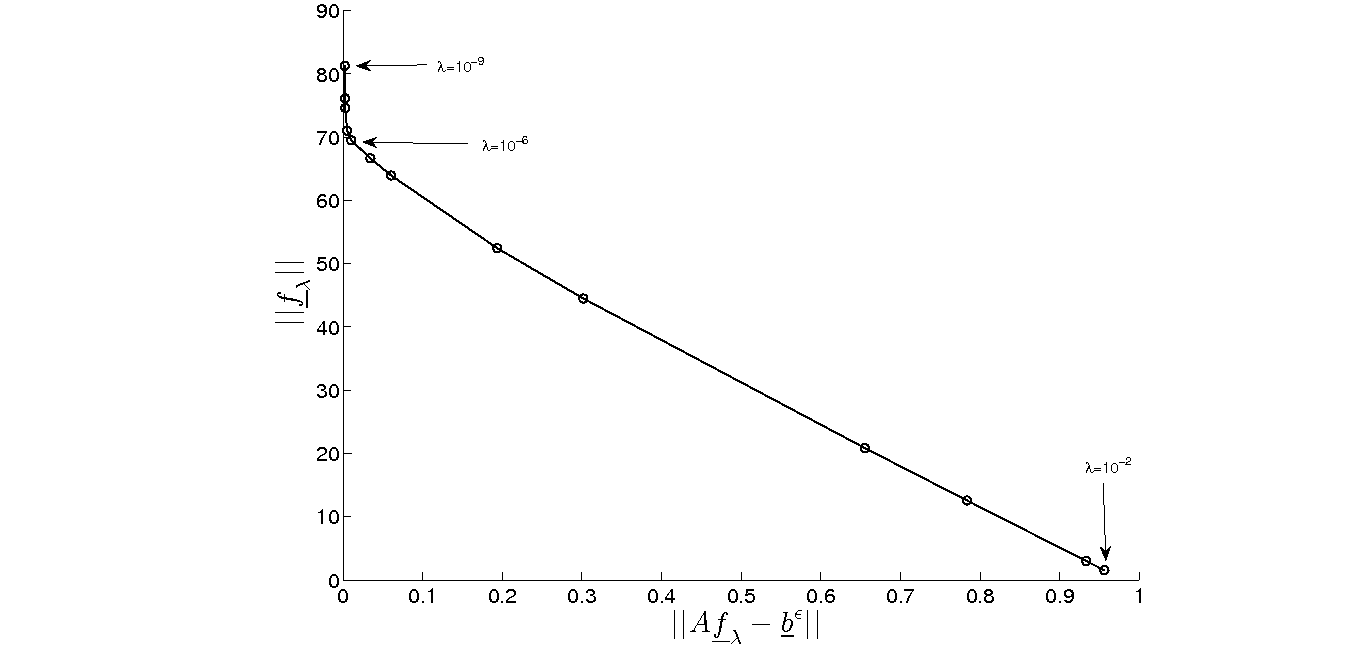}
     \end{center}
\caption{The L-curve for the Tikhonov regularization, for $N=M=80$ and  $p=1\%$ noise, for the inverse problem of Example 1.}
\label{fig:lcurvein0thEx1problem3}
\end{figure}

In order to investigate the stability of the numerical solution we include some ($p=1\%$) noise into the input data \eqref{equx0}, as given by equation \eqref{eq4.5}. The numerical solution for $f(x)$ obtained with $N=M=80$ and no regularization is plotted in Figure \ref{fig:inverseproblemaddnoiseatlambda0Ex1problem3}. It can be clearly seen that very high oscillations appear. This clearly shows that the inverse force problem \eqref{eq12}-\eqref{eq14}, \eqref{eqfsplite} and \eqref{equofxat0t} is ill-posed. In order to deal with this instability we employ the (zeroth-order) Tikhonov regularization which yields the solution 
\begin{eqnarray}
\underline{f}_{\lambda}=(A^{tr}A+\lambda I)^{-1}A^{tr}\underline{b}^{\epsilon},
\label{eqregular} 
\end{eqnarray}
where $I$ is the identity matrix and $\lambda>0$ is a regularization parameter to be prescribed. Including regularization we obtain the numerical solution \eqref{eqregular} whose accuracy error, as a function of $\lambda$, is plotted in Figure \ref{fig:normoffatexactandnumericalEx1problem3}. From this figure it can be 
seen that the minimum of the error occurs around $\lambda=10^{-6}$. Clearly, this
argument cannot be used as a suitable choice for the regularization parameter $\lambda$ in the absence of an analytical (exact) solution \eqref{eq23} being available. However, one possible criterion for choosing $\lambda$ is given by the L-curve method, [6],  
\ \ \quad \\
which plots the residual norm $||A\underline{f}_{\lambda}-\underline{b}^{\epsilon}||$ versus the solution norm $||\underline{f}_{\lambda}||$ for various values of $\lambda$. This is shown in Figure \ref{fig:lcurvein0thEx1problem3} for various values of $\lambda \in \lbrace10^{-9},5\times10^{-9},10^{-8},...,10^{-2}\rbrace.$ The portion to the right of the curve corresponds to large values of $\lambda$ which make the solution oversmooth, whilst the portion to the left of the curve corresponds to small values of $\lambda$ which make the solution undersmooth. The compromise is then achieved around the corner region of the L-curve where the aforementioned portions meet. Figure \ref{fig:lcurvein0thEx1problem3} shows that this corner region includes the values around $\lambda=10^{-6}$, which is a good prediction of the optimal value demonstrated in Figure \ref{fig:normoffatexactandnumericalEx1problem3}.

\begin{figure}[H]
\centering
\includegraphics[width=17cm]{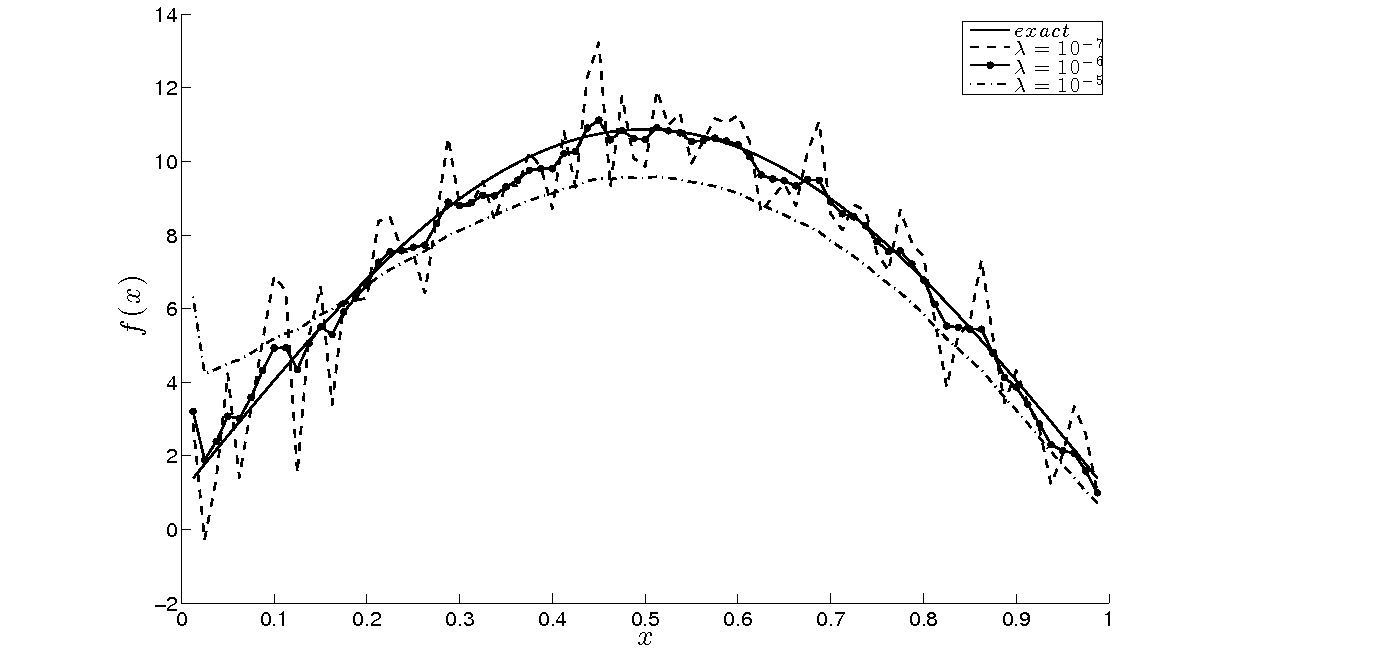}
\caption{The exact solution \eqref{eq23} for the force $f(x)$ in comparison with the numerical solution \eqref{eqregular}, for $N=M=80$, $p=1\%$ noise, and regularization parameters $\lambda \in \lbrace 10^{-7},10^{-6},10^{-5}\rbrace$, for the inverse problem of Example 1.}
\label{fig:optimaloffwithexactEx1problem3}
\end{figure}

Finally, Figure \ref{fig:optimaloffwithexactEx1problem3} shows the regularized numerical solution for $f(x)$ obtained with various values of the regularization parameter $\lambda \in \lbrace 10^{-7},{10}^{-6},10^{-5}\rbrace$ for $p=1\%$ noisy data. From this figure it can be seen that the value of the regularization parameter $\lambda$ can also be chosen by trial and error. By plotting the numerical solution for various values of $\lambda$ we can infer when the instability starts to kick off. For example, in Figure \ref{fig:normoffatexactandnumericalEx1problem3}, the value of $\lambda=10^{-5}$ is too large and the solution is oversmooth, whilst the value of $\lambda=10^{-7}$ is too small and the solution becomes unstable. We could therefore inspect the value of $\lambda=10^{-6}$ and conclude that this is a reasonable choice of the regularization parameter which balances the smoothness with the instability of the solution.
\subsection{Example 2 ($h(x,t)=1+t$)}
This is an example in which we take $h(x,t)=1+t$ a linear function of $t$ and independent of $x$ and consider first the direct problem \eqref{eq12}-\eqref{eq14}  and \eqref{eqfsplite} with the input data
\begin{eqnarray}
u(x,0)=u_{0}(x)=0, \quad
u_{t}(x,0)=v_{0}(x)=0, \quad x\in[0,1],   \label{eq30}
\end{eqnarray}
\begin{eqnarray}
u(0,t)=P_{0}(t)=0,\quad u(1,t)=P_{L}(t)=0, \quad t\in(0,1], \label{eq31}
\end{eqnarray}
\begin{eqnarray}
f(x)=
\begin{cases}
x \ \ \ \ \ \ \ \ \ \ \ \text{if}\ \ \ 0\leq x \leq \frac{1}{2},
\\
1-x \ \ \ \ \ \ \text{if}\ \ \ \frac{1}{2}<x\leq1. \label{eq32}
\end{cases}
\end{eqnarray} 
As in Example 1, since $h(0)=1\neq0$, Theorem 3 ensures the uniqueness of the solution in the class of the functions \eqref{eqth3.2.1}. Also, remark that for this example, the force \eqref{eq32} has a triangular shape, being continuous but non-differentiable at the peak $x=1/2$. This example also does not possess an explicit analytical solution for the displacement $u(x,t)$.

\begin{figure}[H]
    \begin{subfigure}[b]{0.3\textwidth}
     \caption{}
        \centering
        \includegraphics[width=11cm]{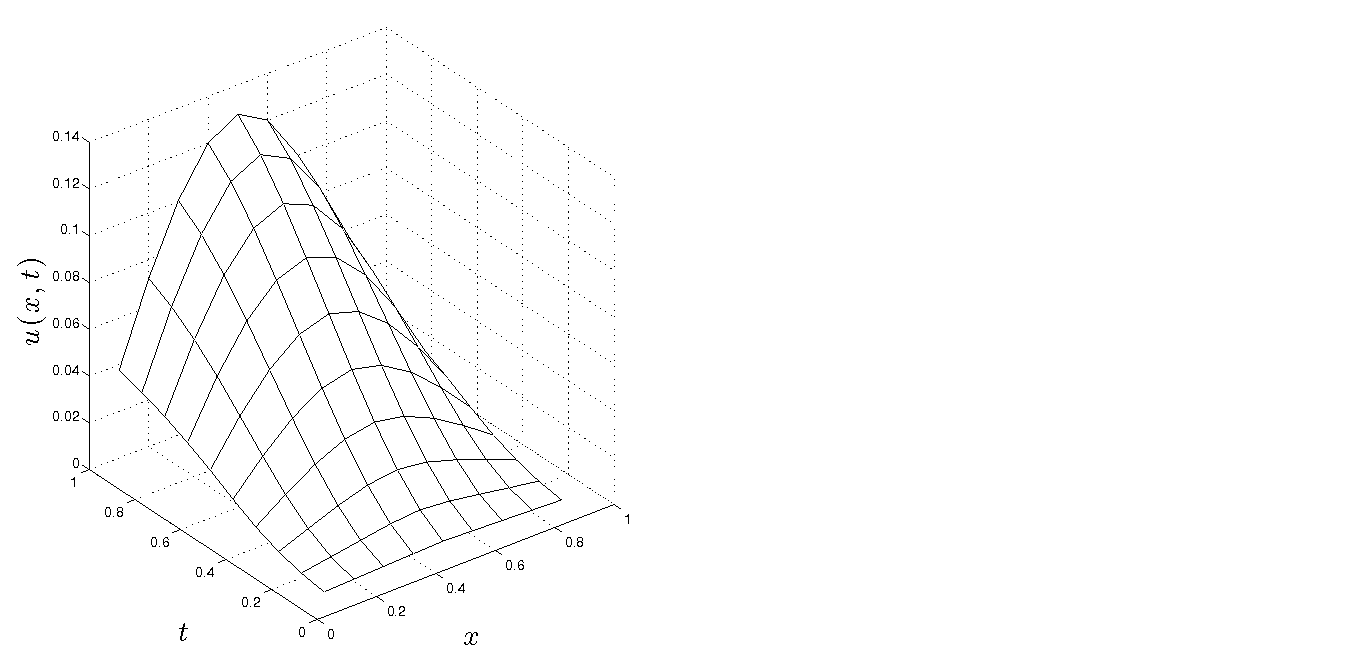}
    \end{subfigure}
    \hspace{2.9cm}
    \begin{subfigure}[b]{0.3\textwidth}
     \caption{}
        \centering
        \includegraphics[width=11cm]{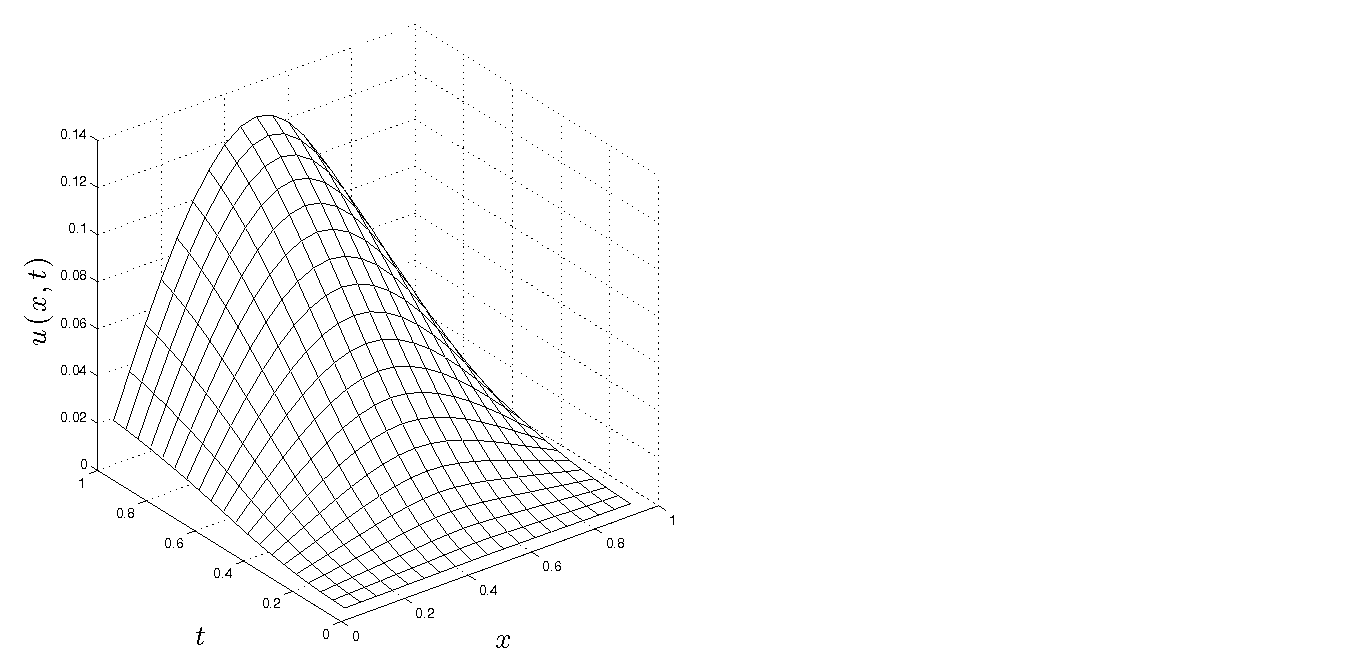}
    \end{subfigure}\\
\\
    \begin{subfigure}[b]{0.3\textwidth}
     \caption{}
        \centering
        \includegraphics[width=11cm]{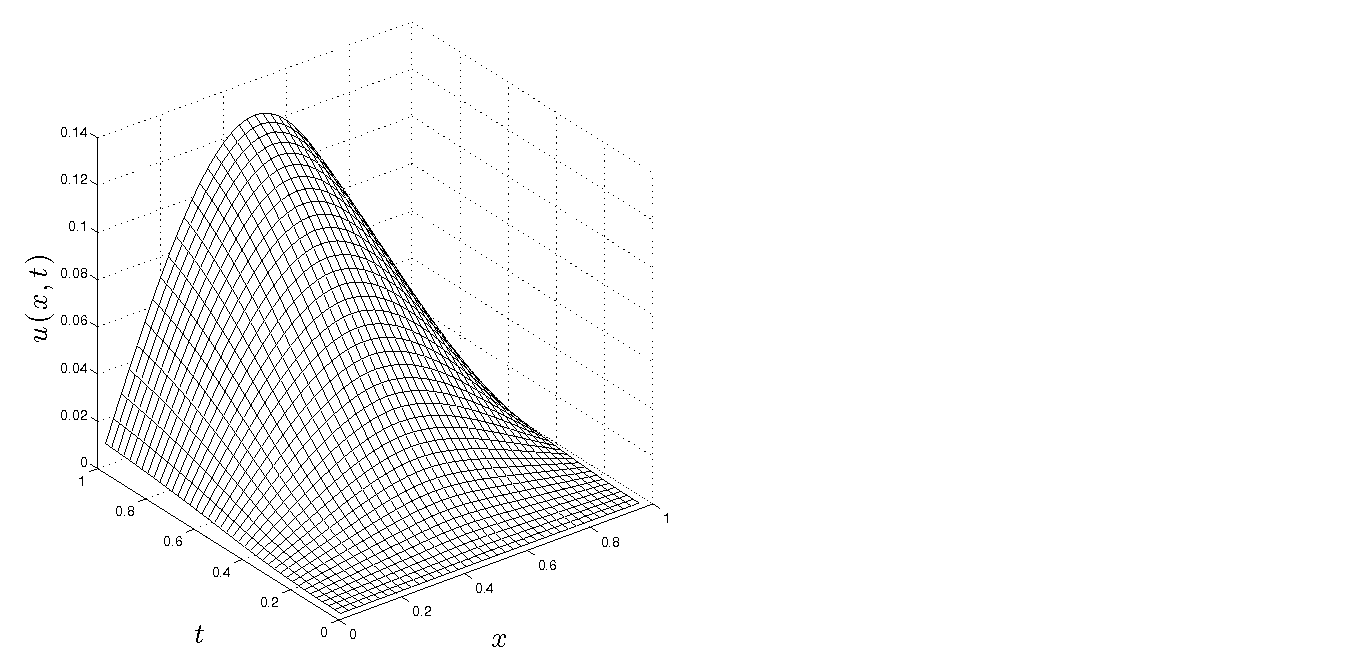}
    \end{subfigure}
\hspace{2.9cm}
    \begin{subfigure}[b]{0.3\textwidth}
     \caption{}
        \centering
        \includegraphics[width=11cm]{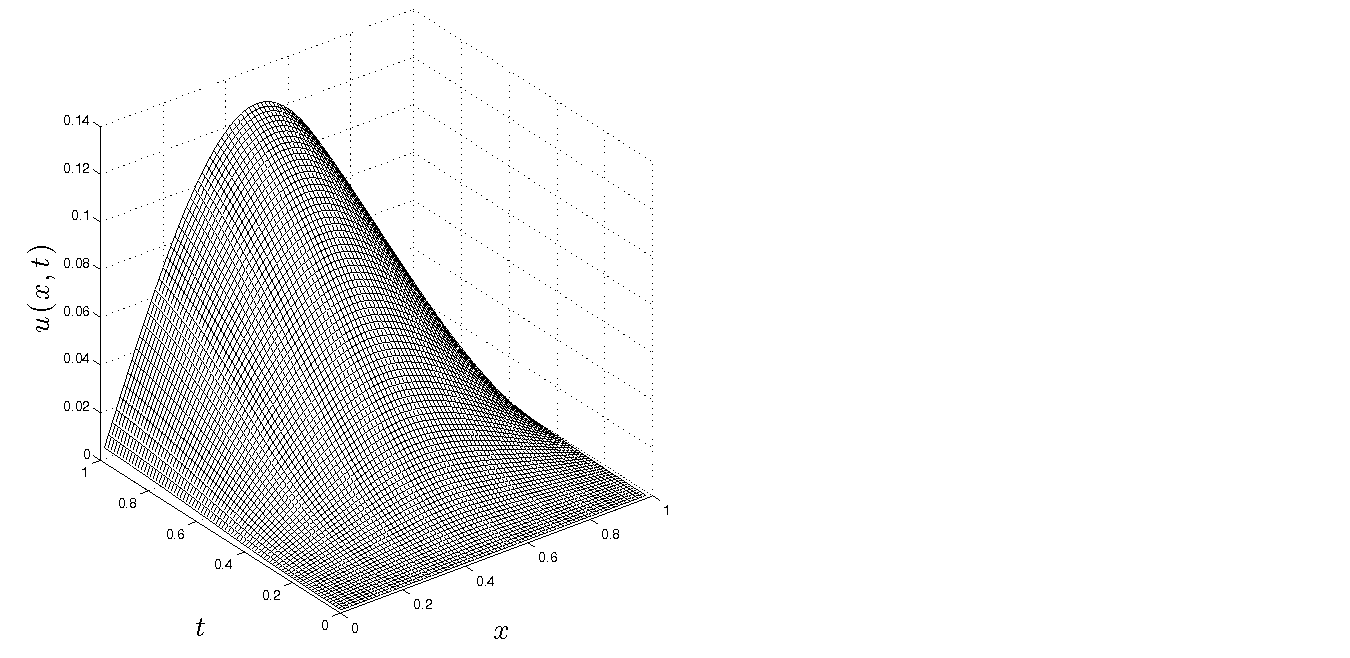}
    \end{subfigure}
    \caption{Numerical solutions for the displacement $u(x,t)$ obtained using the direct problem with various $N=M\in \lbrace 10,20,40,80 \rbrace$ in cases (a)-(d), respectively, for Example 2.}
    \label{fig:directproblemnumericalsolutionofuEx2problem3}
\end{figure}

\begin{figure}[H]
\centering
\includegraphics[width=15cm]{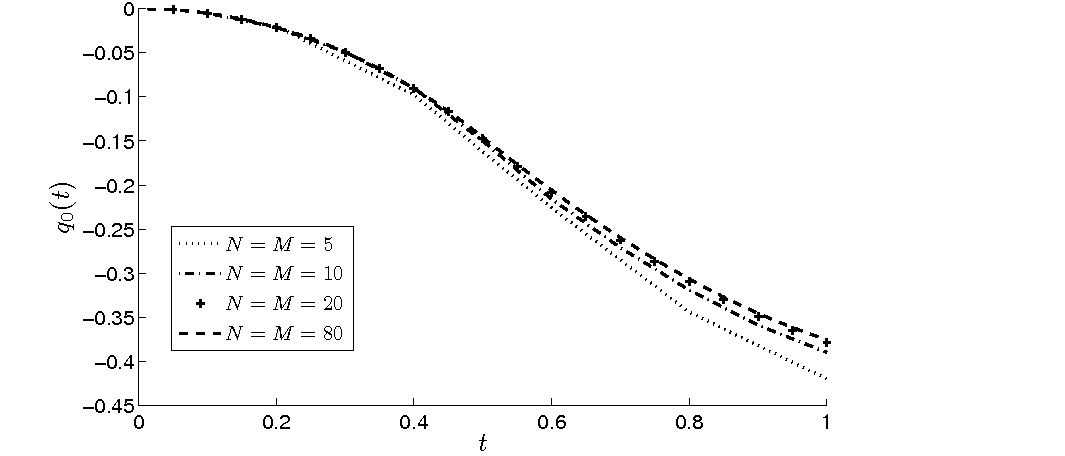}
\caption{Numerical solution for the flux tension at $x=0$, for various $N=M\in\lbrace5,10,20,80\rbrace$, for the direct problem of Example 2.}
\label{fig:numericalsolutionofuxatoEx2problem3}
\end{figure}

The numerical solutions for the displacement $u(x,t)$ at interior points are shown in Figure \ref{fig:directproblemnumericalsolutionofuEx2problem3}. The flux tension \eqref{equofxat0t} is presented in Table 3 and Figure \ref{fig:numericalsolutionofuxatoEx2problem3}. From these figures and table it can be seen that convergent numerical solutions for both $u(x,t)$ and $q_{0}(t)$ are obtained, as $N=M$ increases.
\begin{table}[H]
\caption{The numerical solutions for the flux tension at $x=0$, for the direct problem of Example 2.}
\centering
\begin{tabular}{|c|c|c|c|c|c|c|}
\hline
$t$  &$0.1$   &$0.2$  &$...$   &$0.8$ &$0.9$ &$1$ \\  
\hline
$N=M=10$	&$-0.00500$  &$-0.02100$  &$...$	&$-0.31900$	 &$-0.35900$  &$-0.39000$\\
\hline
$N=M=20$	&$-0.00512$  &$-0.02125$  &$...$	&$-0.3095$	 &$-0.34862$  &$-0.37875$\\
\hline
$N=M=40$	&$-0.00515$  &$-0.02131$  &$...$	&$-0.30712$	 &$-0.34603$  &$-0.37593$\\
\hline
$N=M=80$	&$-0.00516$  &$-0.02132$  &$...$	&$-0.30653$	 &$-0.34538$   &$-0.37523$\\
\hline
\end{tabular}
\end{table}
Consider now the inverse problem given by equations \eqref{eqfsplite} with $h(x,t)=1+t$, equations \eqref{eq30}, \eqref{eq31} and \eqref{equofxat0t} with $q_{0}(t)$ numerically simulated and given in Figure 10 for $N=M=80$. We perturb further this flux by adding to it some $p\in \lbrace1,3,5\rbrace\%$ noise, as given by equation \eqref{eq4.5}. The numerical solution for $f(x)$ obtained with $N=M=80$ and no regularization has been found highly oscillatory and unstable similar to that obtained in Figure \ref{fig:inverseproblemaddnoiseatlambda0Ex1problem3} and therefore is not presented. In order to deal with this instability we employ and test the Tikhonov regularization of various orders such as zero, first and second, which yields the solution, \cite{twos},  
\begin{eqnarray}
\underline{f}_{\lambda}=(A^{tr}A+\lambda D^{tr}_{k}D_{k})^{-1}A^{tr}\underline{b}^{\epsilon},
\label{eqregularex2} 
\end{eqnarray}
where $D_{k}$ is the regularization derivative operator of order $k\in \lbrace0,1,2\rbrace$ and $\lambda  \geq 0$ is the regularization parameter. The regularization derivative operator $D_{k}$ imposes continuity, i.e. class $C^{0}$ for $k=0$, first-order smoothness, i.e. class $C^{1}$ for $k=1$, or second-order smoothness, i.e. class $C^{2}$ for $k=2$. Thus $D_{0}=I$,
\\
\\
$D_{1}=
\begin{pmatrix}
1 & -1 & 0  &0 & ... & 0  \\
0 & 1 & -1 & 0& ... & 0 \\
... & ... & ... & ... & ...&...\\
0 & 0 & ... & 0 & 1 &-1
\end{pmatrix}$, \ $D_{2}=
\begin{pmatrix}
1 & -2 & 1 & 0 & 0  & ... & 0  \\
0 & 1 & -2 & 1 & 0 &... & 0 \\
... & ... & ... & ... & ...& ... & ...&\\
0 & 0 & ... & 0 & 1 & -2 & 1
\end{pmatrix}$.
\\
\\
Observe that for $k=0$, equation \eqref{eqregularex2} becomes the zeroth-order regularized solution \eqref{eqregular} which was previously employed in Example 1 in order to obtain a stable solution.

\begin{figure}[H]
    \begin{subfigure}[b]{0.3\textwidth}
       \caption{}
        \includegraphics[width=11cm]{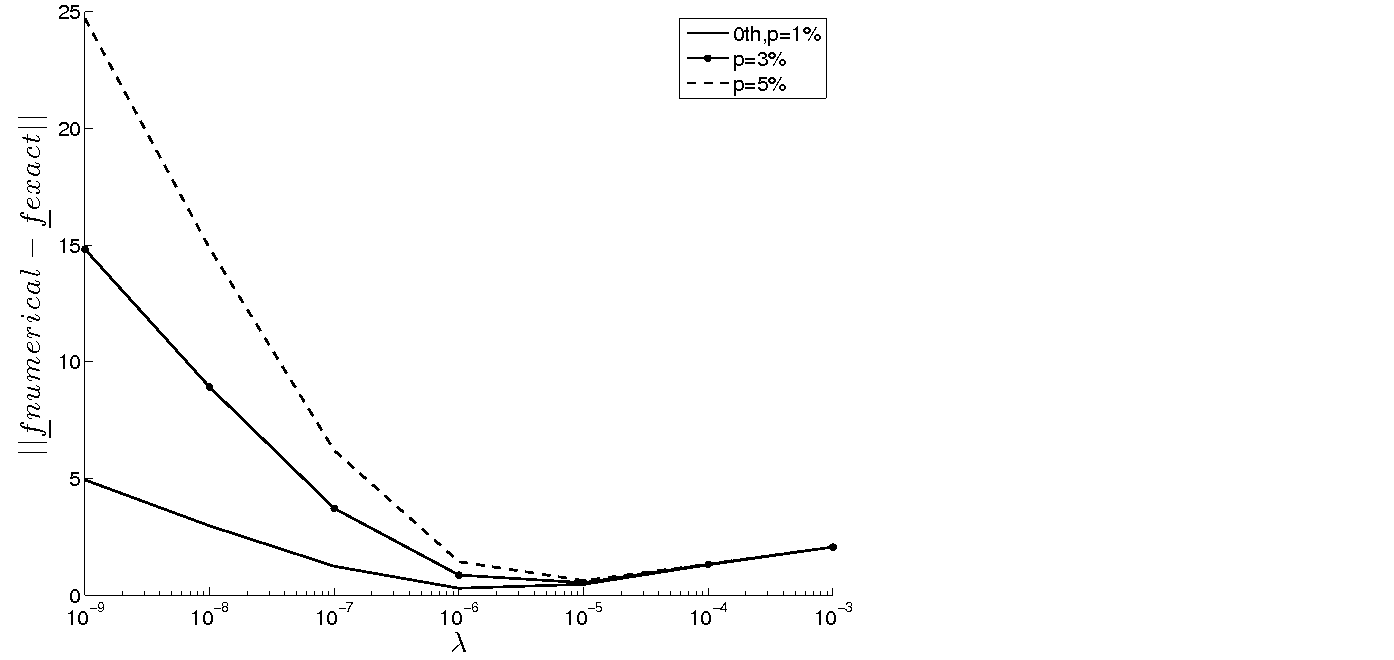}
    \end{subfigure}
    \hspace{3cm}
    \begin{subfigure}[b]{0.3\textwidth}
     \caption{\ \ \ \ \ }
        \includegraphics[width=11cm]{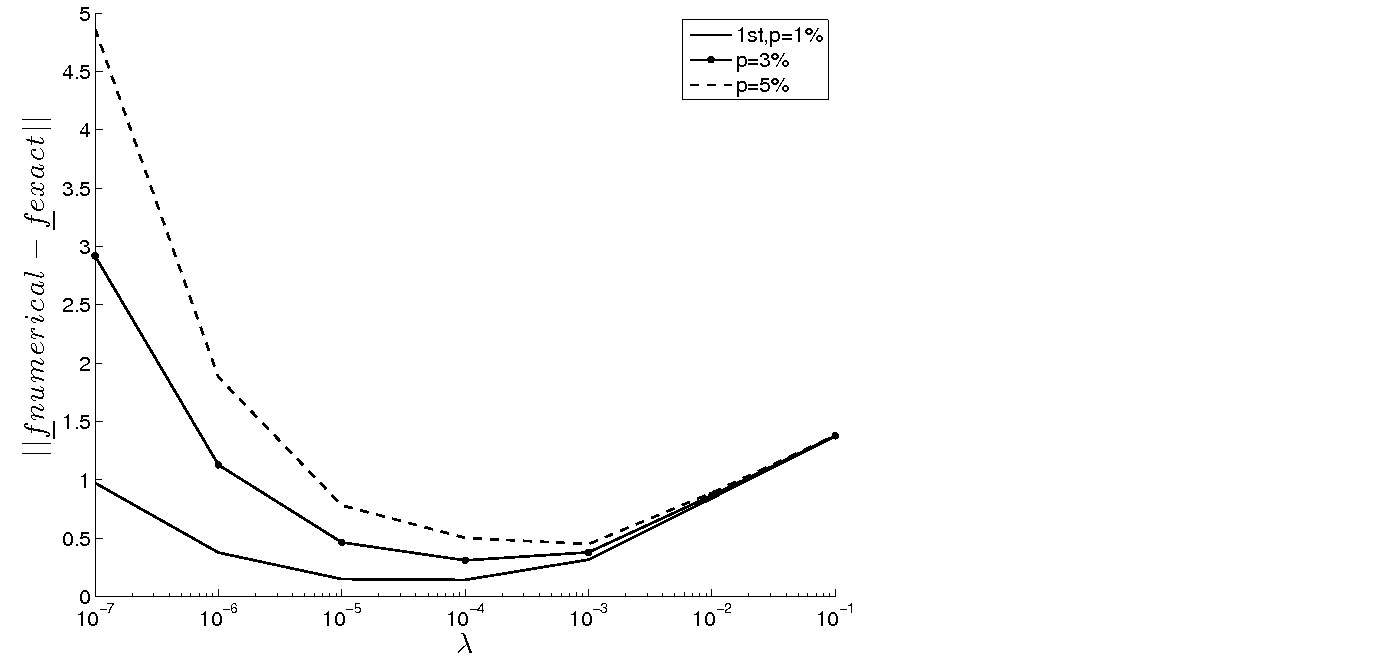}
    \end{subfigure}\\
   \begin{center}
    \begin{subfigure}[b]{0.3\textwidth}
     \caption{}
        \includegraphics[width=11cm]{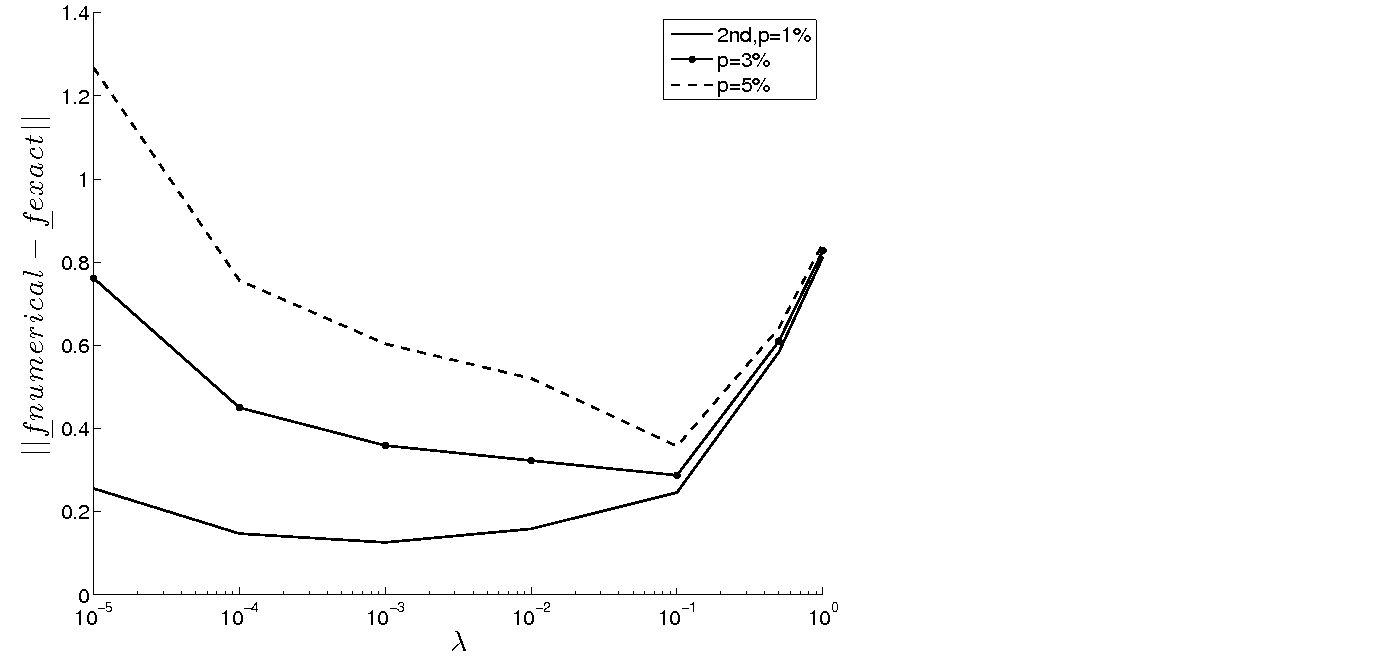}
    \end{subfigure}
\end{center}
   \caption{The accuracy error $||\underline{f}numerical-\underline{f}exact||$, as a function of $\lambda$, for $M=N=80$, $p\in \lbrace1,3,5\rbrace\%$ noise, obtained using (a) zeroth, (b) first, and (c) second-order regularization, for the inverse problem of Example 2.}
   \label{fig:optimaloffwithexactone0th1st2ndEx2problem3}
\end{figure}

\begin{figure}[H]
        \centering
      \includegraphics[width=17cm]{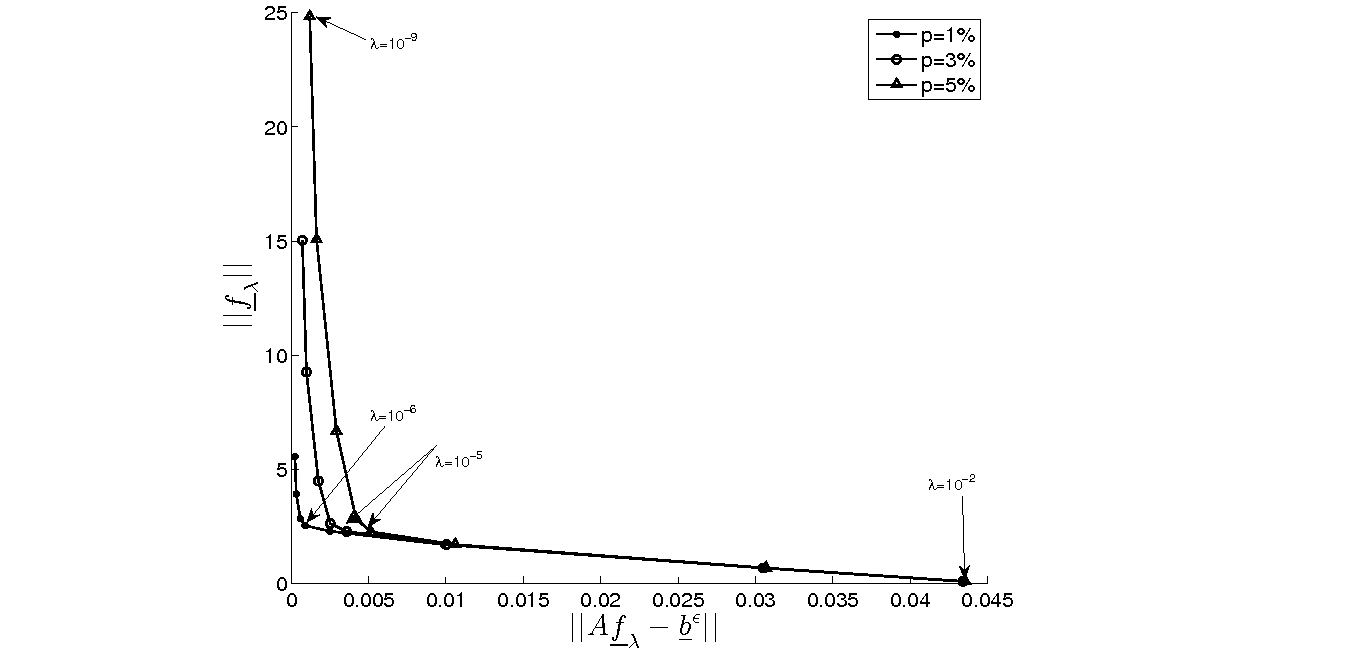}
\caption{The L-curve for the zeroth-order Tikhonov regularization, for $N=M=80$ and  $p\in \lbrace1,3,5\rbrace\%$ noise, for  the inverse problem of Example 2.}
\label{fig:lcurve0thEx2problem3}
\end{figure}

Including regularization we obtain the solution \eqref{eqregularex2} whose accuracy error, as a function of $\lambda$, is plotted in Figure \ref{fig:optimaloffwithexactone0th1st2ndEx2problem3} for various orders of regularization $k\in \lbrace 0,1,2 \rbrace$. From this figure it can be seen that there are wide ranges for choosing the regularization parameters in the valleys of minima of the plotted error curves. The minimum points $\lambda_{opt}$ and the corresponding accuracy errors are listed in Table 4. The L-curve criterion for choosing $\lambda$ in the zeroth-order regularisation is shown in Figure \ref{fig:lcurve0thEx2problem3} for various values of $\lambda \in \lbrace 10^{-9},10^{-8},...,10^{-2}\rbrace$ and for $p\in \lbrace1,3,5\rbrace\%$ noisy data. This figure shows that the L-corner
region includes the values around $\lambda=10^{-6}$ for $p=1\%$, $\lambda=10^{-5}$ for $p=3\%$, and $\lambda=10^{-5}$ for $p=5\%$. Similar L-curves, which plot the penalised solution norm $||D_{k}\underline{f}_{\lambda}||$ versus the residual norm $||A\underline{f}_{\lambda}-\underline{b}^{\epsilon}||$, have been obtained for the first and second-order regularizations and therefore they are not illustrated.   

\begin{figure}[H]
    \begin{subfigure}[b]{0.3\textwidth}
       \caption{}
        \includegraphics[width=11cm]{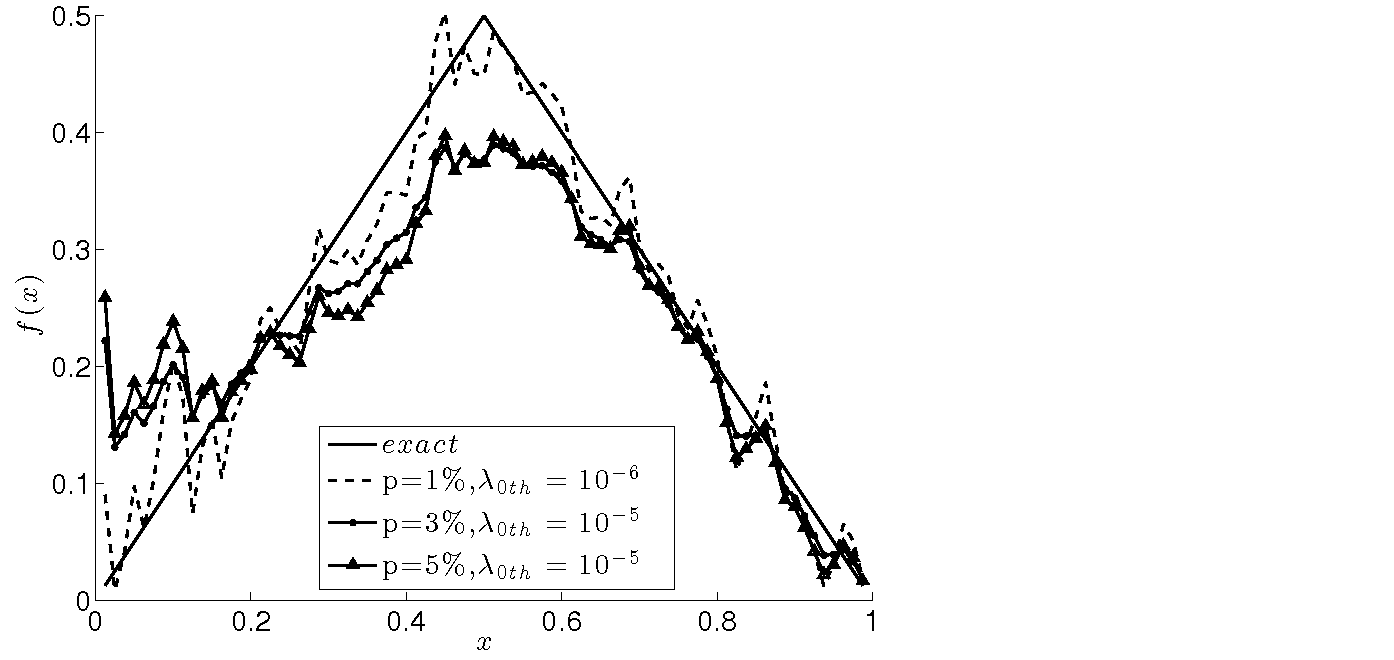}
    \end{subfigure}
    \hspace{3.9cm}
    \begin{subfigure}[b]{0.3\textwidth}
        \caption{}
        \includegraphics[width=11cm]{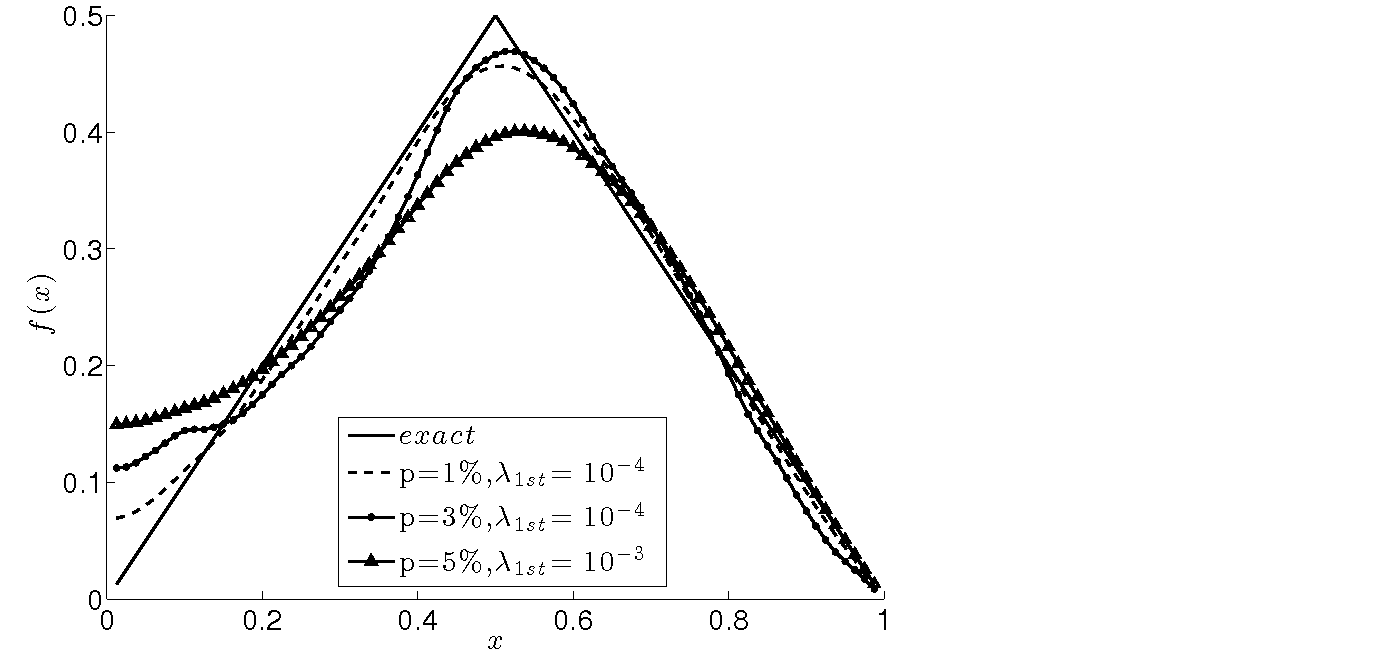}
    \end{subfigure}\\
\begin{center}
    \begin{subfigure}[b]{0.3\textwidth}
     \caption{}
        \includegraphics[width=11cm]{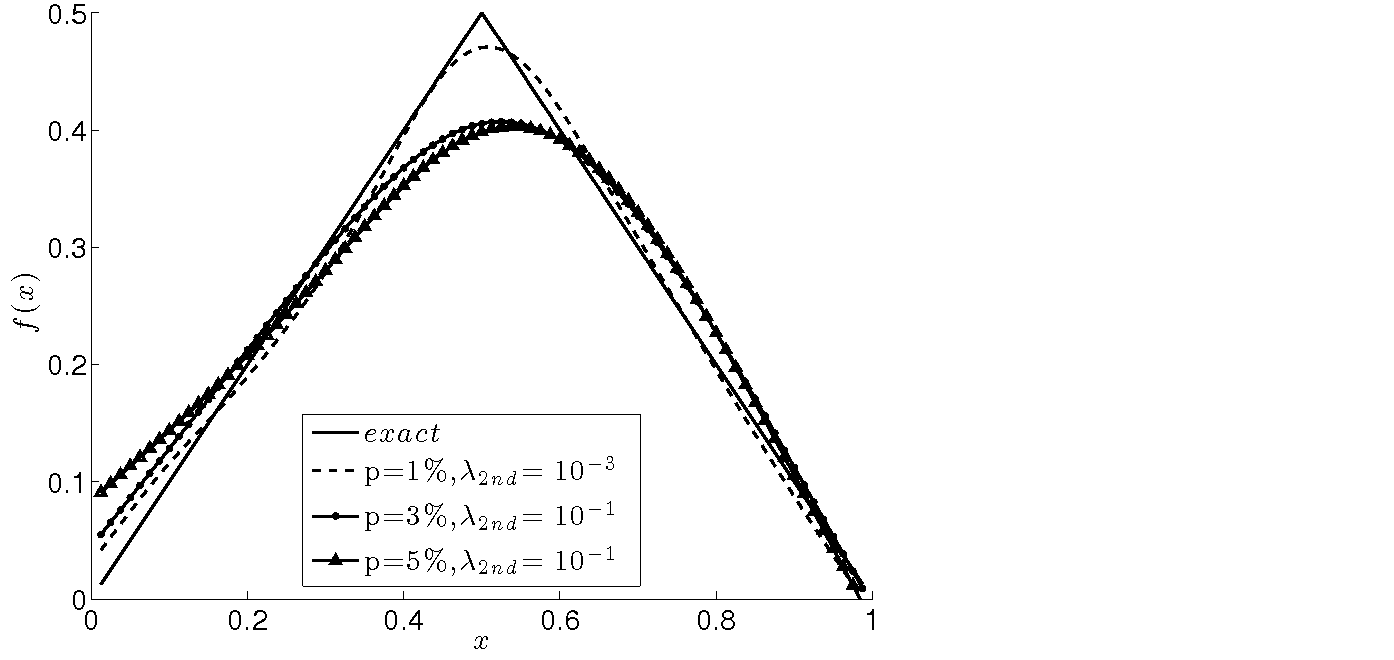}
    \end{subfigure}
\end{center}
\caption{The exact solution \eqref{eq32} for the force $f(x)$ in comparison with the numerical regularized solution \eqref{eqregularex2}, for $N=M=80$, $p\in\lbrace1,3,5\rbrace\%$ noise, and various order regularization methods, for the inverse problem of Example 2.}
\label{fig:optimaloffwithexact0th1st2ndEx2problem3}
\end{figure}

Figure \ref{fig:optimaloffwithexact0th1st2ndEx2problem3} shows the regularized numerical solutions \eqref{eqregularex2} for $f(x)$ obtained with the values of the regularization parameter $\lambda_{opt}$ given in Table 4 for $p\in\lbrace1,3,5\rbrace\%$ noisy data.  From this figure it can be seen that the numerical results are stable and they become more accurate as the amount of noise $p$ decreases.  
\begin{table}[H]
\caption{The accuracy error $||\underline{f}numerical-\underline{f}exact||$ for various order regularization methods and percentages of noise $p$, for the inverse problem of Example 2. The values of $\lambda_{opt}$ are also included.}
\centering
\begin{tabular}{|c|c|c|c|c|c|c|}
\hline
Regularization  &$p=1\%$   &$p=3\%$  &$p=5\%$ \\  
\hline
zeroth	&$\lambda_{opt}=10^{-6}$  &$\lambda_{opt}=10^{-5}$  &$\lambda_{opt}=10^{-5}$\\
$$ & $0.2987$ & $0.5389$ & $0.6259$ \\
\hline
first	&$\lambda_{opt}=10^{-4}$  &$\lambda_{opt}=10^{-4}$  &$\lambda_{opt}=10^{-3}$\\
$$ & $0.1433$ & $0.3112$ & $0.4494$ \\
\hline
second 	&$\lambda_{opt}=10^{-3}$  &$\lambda_{opt}=10^{-1}$  &$\lambda_{opt}=10^{-1}$\\
$$ & $0.1264$ & $0.2876$ & $0.3576$ \\
\hline
\end{tabular}
\end{table}  
\subsection{Example 3 ($h(x,t)=1+x+t$)}
All the data and details of the numerical implementation are the same as those for Example 2, except that for the present example $h(x,t)=1+x+t$ in equation \eqref{eqfsplite}. Since in this case $h$ depends also on $x$ we cannot apply Theorem 3, but we can apply instead Theorem 1, because $H=0$ in \eqref{eqth1.1} is sufficiently small. This then ensures the uniqueness of the solution in the class of functions \eqref{eqth1.2}, which in $n=1$-dimension reads as
\begin{eqnarray}
u\in L^{2}(0,T;H^{1}(0,L)), \quad u_{t}\in L^{2}(0,T;L^{2}(0,L)), \quad u_{tt}\in L^{2}(0,T;(H^{1}(0,L))^{\prime}), \notag \\
f\in L^{2}(0,L). \quad \quad \qquad \quad \quad \qquad \quad \quad \qquad \quad \quad \qquad \quad \quad \qquad \quad \quad \qquad \quad \quad \qquad \label{eqth1.2.1}
\end{eqnarray}

\begin{figure}[H]
    \begin{subfigure}[b]{0.3\textwidth}
       \caption{}
        \includegraphics[width=11cm]{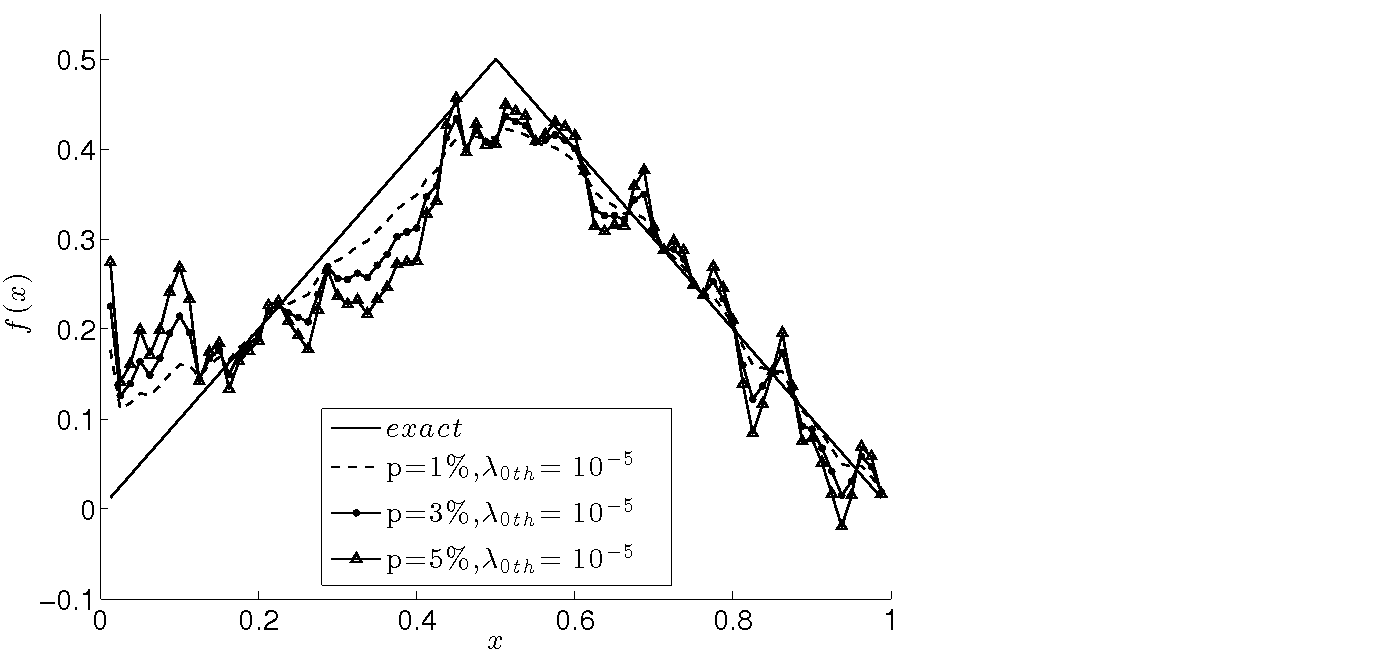}
    \end{subfigure}
    \hspace{3.9cm}
    \begin{subfigure}[b]{0.3\textwidth}
        \caption{}
        \includegraphics[width=11cm]{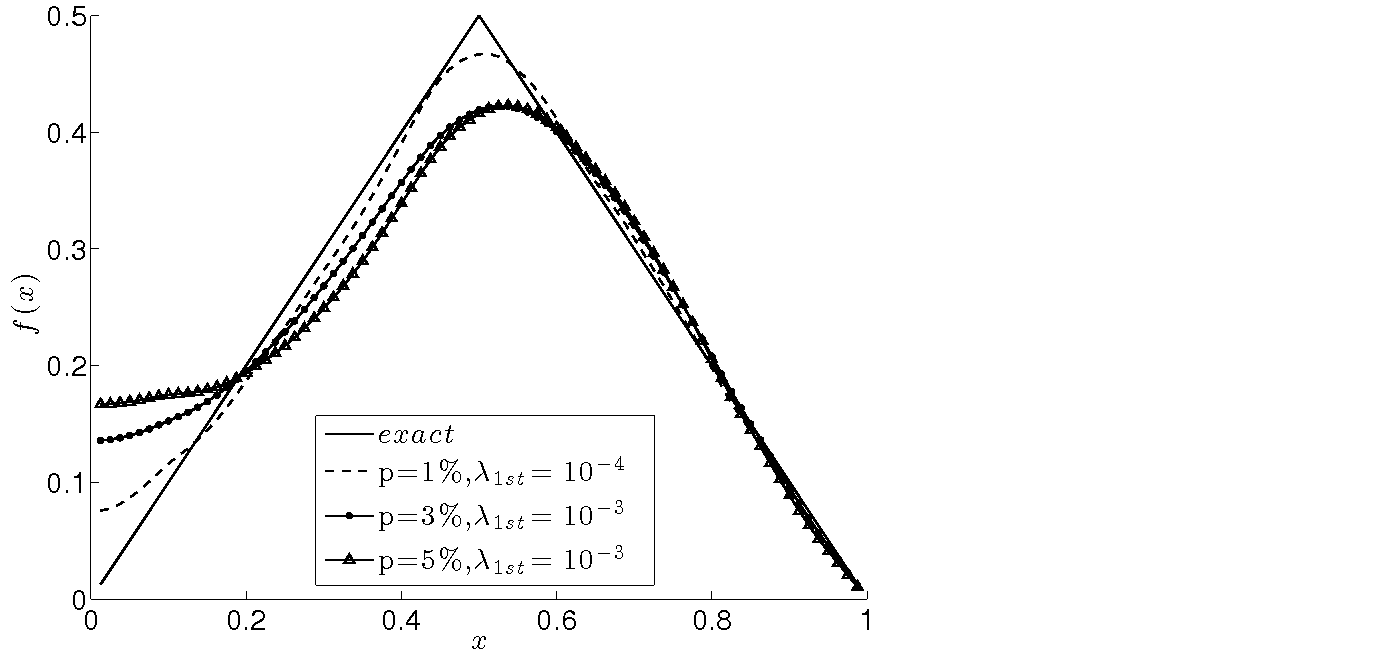}
    \end{subfigure}\\
\begin{center}
    \begin{subfigure}[b]{0.3\textwidth}
     \caption{}
        \includegraphics[width=11cm]{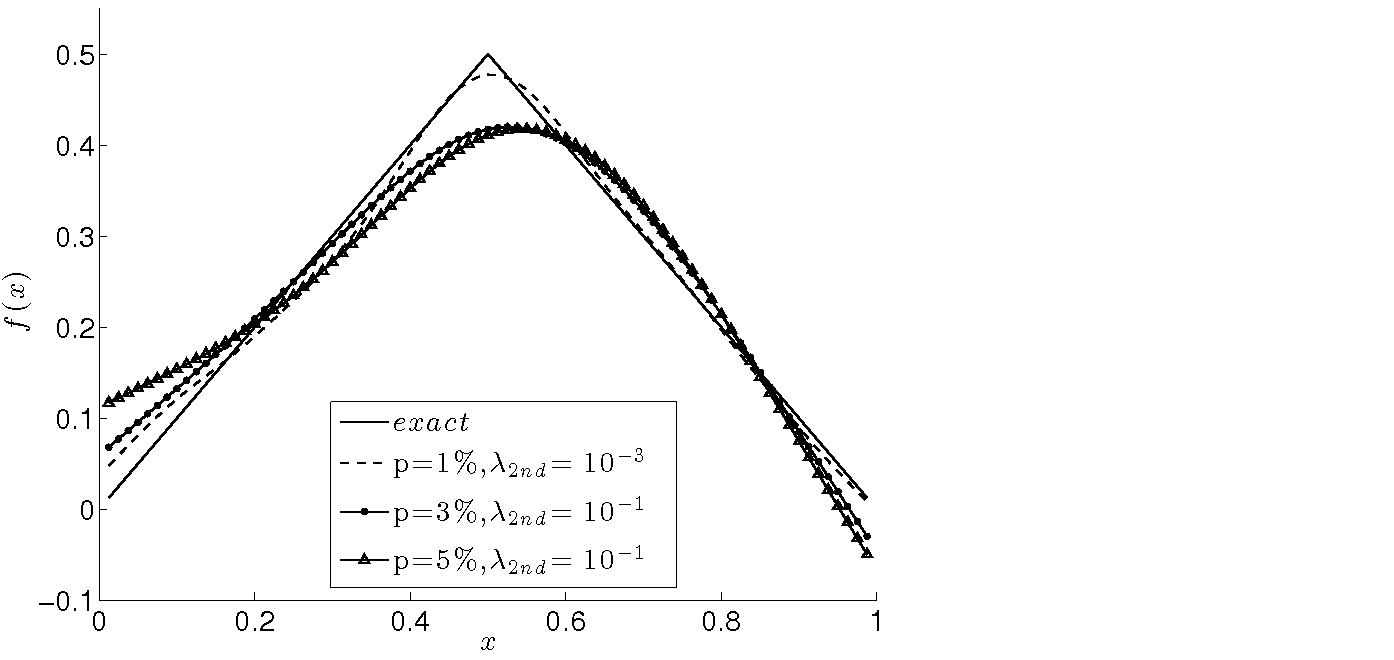}
    \end{subfigure}
\end{center}
\caption{The exact solution \eqref{eq32} for the force $f(x)$ in comparison with the regularized numerical solution \eqref{eqregularex2}, for $N=M=80$, $p\in\lbrace1,3,5\rbrace\%$ noise, and various order regularization methods, for the inverse problem of Example 3.}
\label{fig:inverseoptimalfwithexact0th1st2ndEx3problem3}
\end{figure}

Figure \ref{fig:inverseoptimalfwithexact0th1st2ndEx3problem3} shows the regularized numerical solution for $f(x)$ obtained with various values of the regularization parameters listed in Table 5 for $p\in\lbrace1,3,5\rbrace\%$ noisy data. From this figure it can be seen that stable numerical solutions are obtained.
\begin{table}[H]
\caption{The accuracy error $||\underline{f}numerical-\underline{f}exact||$ for various order regularization methods and percentages of noise $p$, for the inverse problem of Example 3. The values of $\lambda_{opt}$ are also included.}
\centering
\begin{tabular}{|c|c|c|c|c|c|c|}
\hline
Regularization  &$p=1\%$   &$p=3\%$  &$p=5\%$ \\  
\hline
zeroth	&$\lambda_{opt}=10^{-5}$  &$\lambda_{opt}=10^{-5}$  &$\lambda_{opt}=10^{-5}$\\
$$ & $0.35490$ & $0.49093$ & $0.65283$ \\
\hline
first	&$\lambda_{opt}=10^{-4}$  &$\lambda_{opt}=10^{-3}$  &$\lambda_{opt}=10^{-3}$\\
$$ & $0.14821$ & $0.35679$ & $0.45932$ \\
\hline
second 	&$\lambda_{opt}=10^{-3}$  &$\lambda_{opt}=10^{-1}$  &$\lambda_{opt}=10^{-1}$\\
$$ & $0.13326$ & $0.27424$ & $0.39021$ \\
\hline
\end{tabular}
\end{table}  
\subsection{Example 4 ($h(x,t)=t^{2}$)}
All the details are the same as those for Example 2, except that for the present example $h(x,t)=t^{2}$ in equation \eqref{eqfsplite} is independent of $x$, but is a nonlinear function of $t$. Furthermore, one can see that $h(0)=0$ and also, condition \eqref{eqth1.1} is violated. Hence, we cannot apply the uniqueness Theorems 1-3 and, in this case, we expect a more severe situation than in the previous examples to occur. This is reflected in the very large condition numbers of the matrix $A$ reported in Table 1 for Example 4 in comparison with the milder condition numbers for Examples 1-3. 

\begin{figure}[H]
\includegraphics[width=19cm]{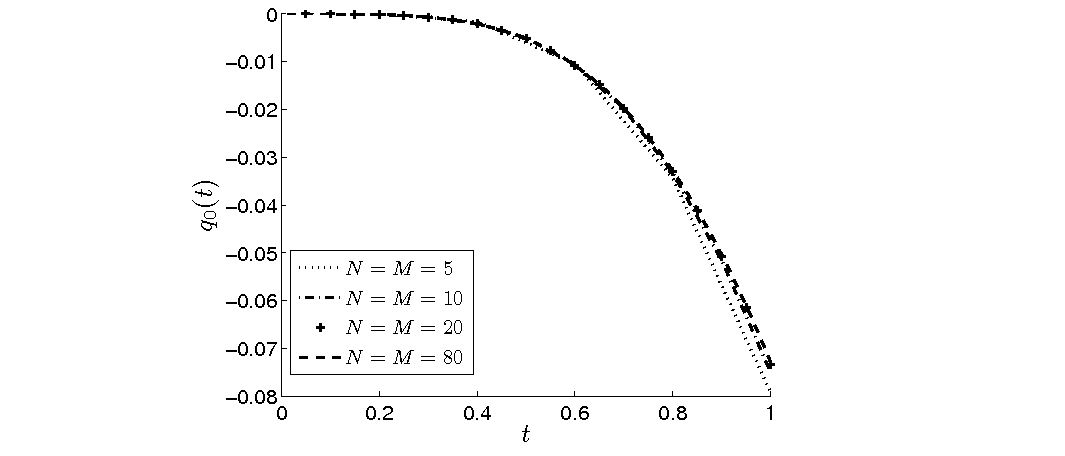}\\
\caption{Numerical solution for the flux tension at $x=0$, for various $N=M\in\lbrace5,10,20,80\rbrace$, for the direct problem of Example 4.}
\label{fig:directproblemnumericalsolutionuxat0Ex4problem3}
\end{figure}

The numerical solution for the flux tension \eqref{equofxat0t} obtained by solving the direct problem given by equation \eqref{eqfsplite} with $h(x,t)=t^{2}$ and equations \eqref{eq30}-\eqref{eq32} is illustrated in Figure \ref{fig:directproblemnumericalsolutionuxat0Ex4problem3} for various mesh sizes. From this figure it can be seen that a rapidly convergent numerical solution is achieved. As in Example 2, we add noise to the numerical flux $q_{0}(t)$ obtained with the finer mesh $N=M=80$.

\begin{figure}[H]
    \begin{subfigure}[b]{0.3\textwidth}
       \caption{}
        \includegraphics[width=11cm]{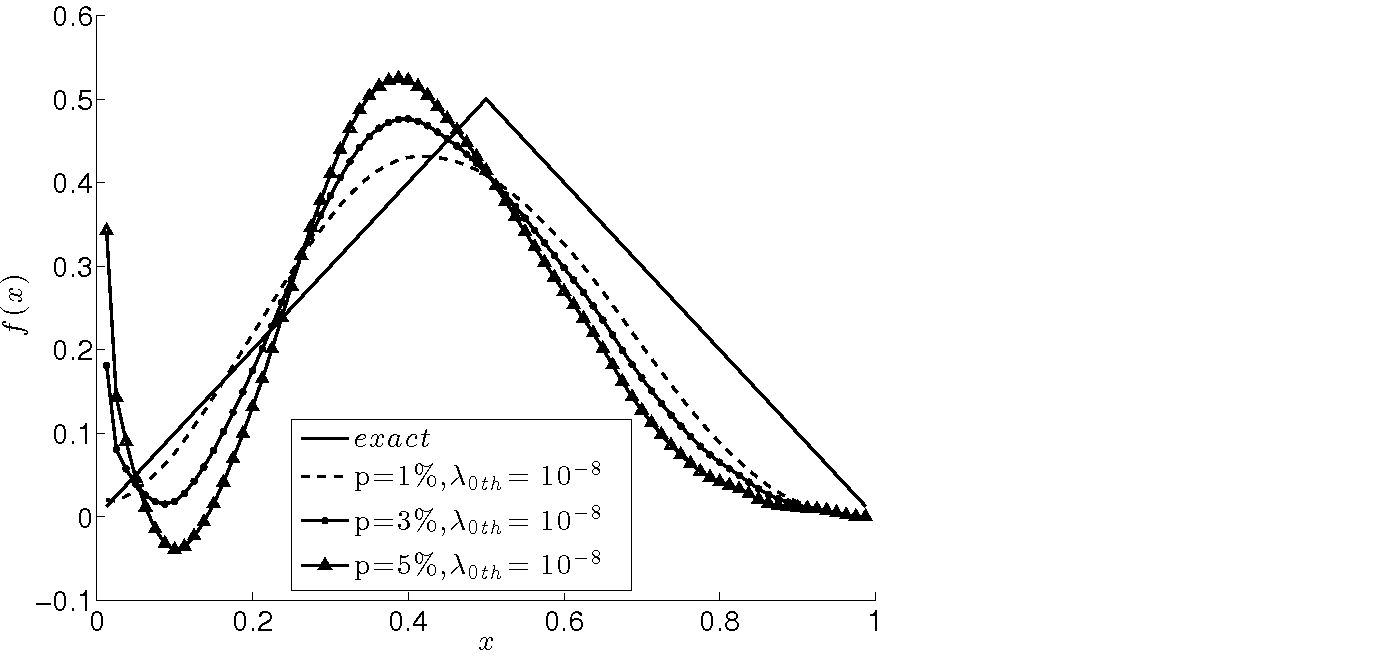}
    \end{subfigure}
    \hspace{3.9cm}
    \begin{subfigure}[b]{0.3\textwidth}
        \caption{}
        \includegraphics[width=11cm]{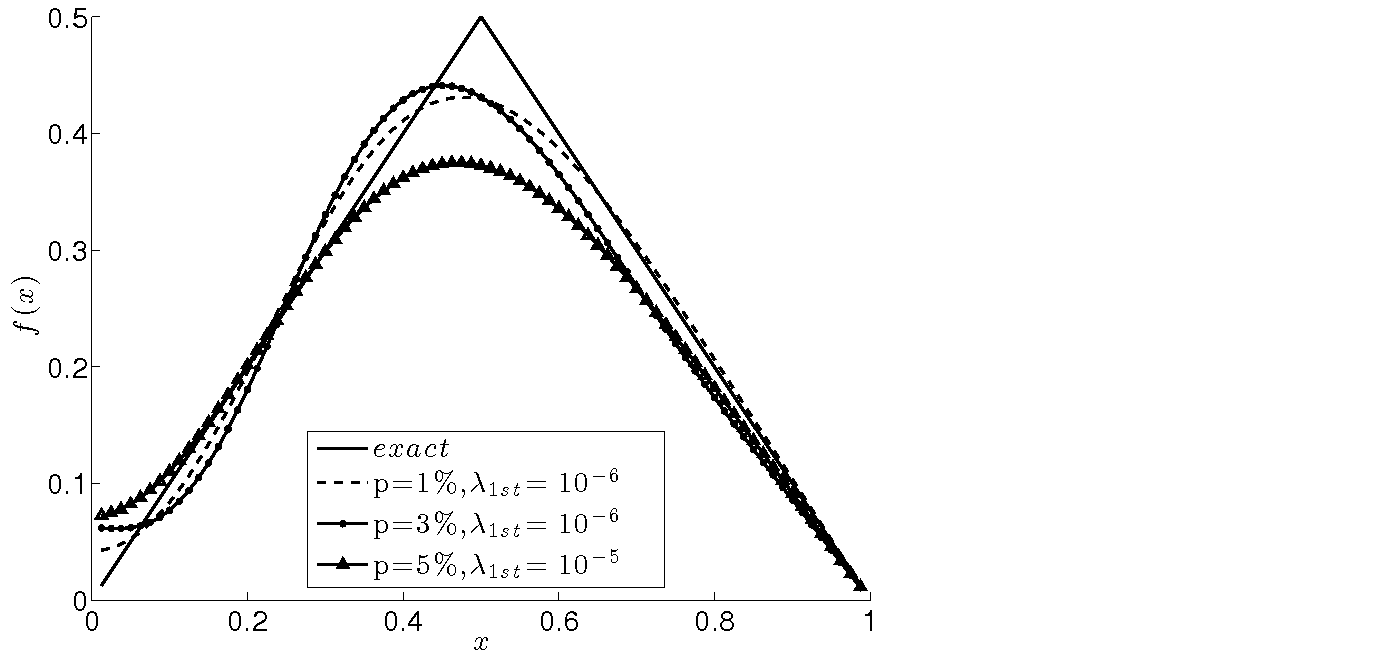}
    \end{subfigure}\\
\begin{center}
    \begin{subfigure}[b]{0.3\textwidth}
     \caption{}
        \includegraphics[width=11cm]{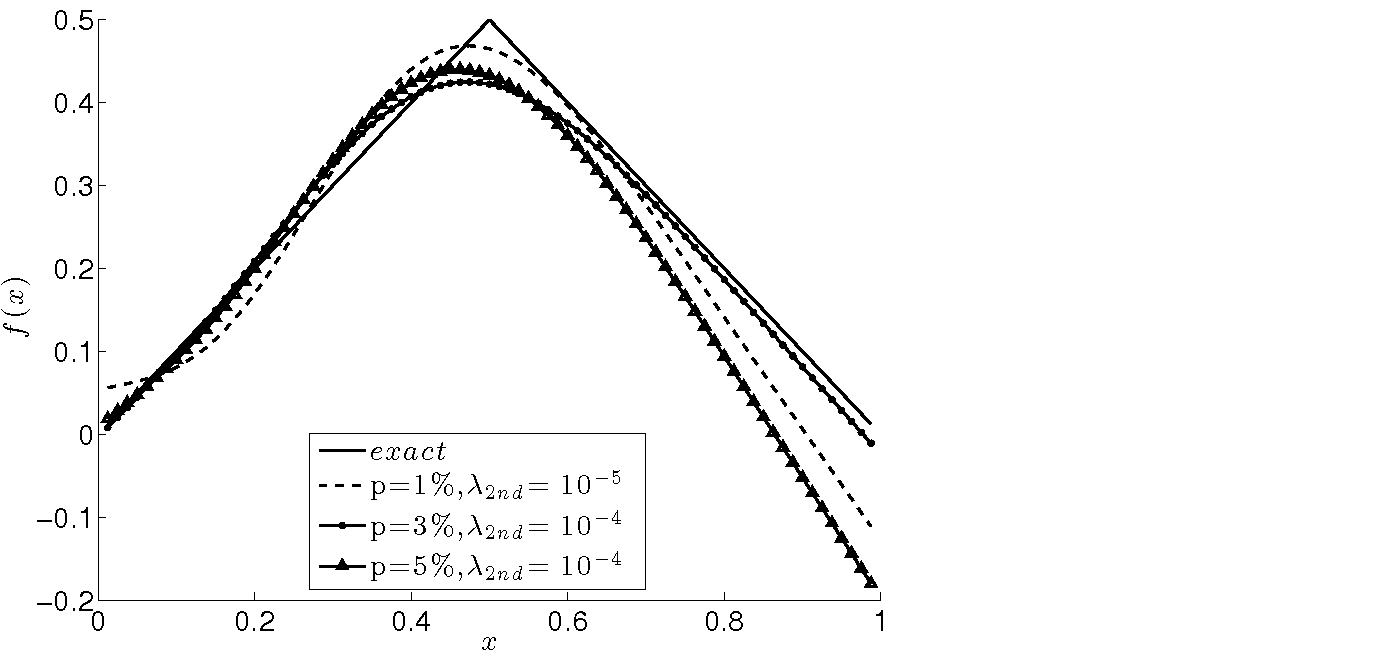}
    \end{subfigure}
\end{center}
\caption{The exact solution \eqref{eq32} for the force $f(x)$ in comparison with the regularized numerical solution \eqref{eqregularex2}, for $N=M=80$ and $p\in\lbrace1,3,5\rbrace\%$ noise, and various order regularization methods, for the inverse problem of Example 4.}
\label{fig:optimaloffwithexact0th1st2ndEx4problem3}
\end{figure}

Figure \ref{fig:optimaloffwithexact0th1st2ndEx4problem3} shows the regularized numerical solution for $f(x)$ obtained with various regularization parameters listed in Table 6 for
$p\in \lbrace1,3,5\rbrace\%$ noisy data. As in all the previous examples, stable numerical solutions are obtained. However, in contrast to Examples 2 and 3, the first-order regularization seems to perform better than the second-order regularization, with the latter one also presenting some unexpected behaviour of increase in accuracy when $p$ increases from $1\%$ to $3\%$. These conclusions may be attributed  to the severe ill-posedness of the Example 4 which, as discussed above, in addition to ill-conditioning it fails to satisfy the conditions for uniqueness of solution of Theorems 1-3.
\begin{table}[H]
\caption{The accuracy error $||\underline{f}numerical-\underline{f}exact||$ for various order regularization methods and percentages of noise $p$, for the inverse problem of Example 4. The values of $\lambda_{opt}$ are also included.}
\centering
\begin{tabular}{|c|c|c|c|c|c|c|}
\hline
Regularization  &$p=1\%$   &$p=3\%$  &$p=5\%$ \\  
\hline
zeroth	&$\lambda_{opt}=10^{-8}$  &$\lambda_{opt}=10^{-8}$  &$\lambda_{opt}=10^{-8}$\\
$$ & $0.5947$ & $0.8082$ & $1.0863$ \\
\hline
first	&$\lambda_{opt}=10^{-6}$  &$\lambda_{opt}=10^{-6}$  &$\lambda_{opt}=10^{-5}$\\
$$ & $0.1826$ & $0.2668$ & $0.4053$ \\
\hline
second 	&$\lambda_{opt}=10^{-5}$  &$\lambda_{opt}=10^{-4}$  &$\lambda_{opt}=10^{-4}$\\
$$ & $0.4313$ & $0.2178$ & $0.6912$ \\
\hline
\end{tabular}
\end{table}  
\section{Extension to Multiple Sources}
In this section, we consider an extension of the inverse space-dependent problem, in the situation when 
\begin{eqnarray}
F(\underline{x},t)=f(\underline{x})h(\underline{x},t)+g(\underline{x})\theta(\underline{x},t), \quad (\underline{x},t)\in\Omega\times(0,T]. \label{eqextension1}
\end{eqnarray}
where $h(\underline{x},t)$ and $\theta(\underline{x},t)$ are given functions and $f(\underline{x})$ and $g(\underline{x})$ are space-dependent unknown force components to be determined. Under the assumption \eqref{eqextension1}, equation \eqref{eq1} in one-dimension, i.e. $n=1$ and $\Omega=(0,L)$, becomes
\begin{eqnarray}
u_{tt}(x,t)=u_{xx}(x,t)+f(x)h(x,t)+g(x)\theta(x,t), \quad (x,t)\in(0,L)\times(0,T]. \label{eqfandgsplite}
\end{eqnarray}
This has to be solved subject to the initial and boundary conditions \eqref{eq12}-\eqref{eq14} and the overspecified flux tensions at both ends of the string, namely, \eqref{equofxat0t} and
 \begin{eqnarray}
 \frac{\partial{u}}{\partial{x}}(L,t)=q(L,t)=:q_{L}(t), \quad t\in(0,T]. \label{equofxatLt}
 \end{eqnarray}
Then uniqueness of solution still holds in the case $h(x,t)=1$, $\theta(x,t)=t$, see Theorem 8 of \cite{engl94}, but for more general cases, e.g. $h(x,t)=1$, $\theta(x,t)=t^{2}$, the solution ($f(x),g(x),u(x,t)$) is not unique, see the counterexample to uniqueness given in \cite{engl94}. 
 
In discretised finite-difference form equations \eqref{eq12}-\eqref{eq14} and \eqref{eqfandgsplite}  recast as equations \eqref{eq16.1}, \eqref{eq16.2},
\begin{eqnarray}
u_{i,j+1}-(\delta t)^{2}f_{i}h_{i,j}-(\delta t)^{2}g_{i}\theta_{i,j}=r^{2}u_{i+1,j}+2(1-r^{2})u_{i,j}+r^{2}u_{i-1,j}-u_{i,j-1}, \label{eqfandgdmsplitef}\\
\quad \quad  \quad i=\overline{1,(M-1)}, \quad j=\overline{1,(N-1)},\notag 
\end{eqnarray}
and
\begin{eqnarray}
&&u_{i,1}-\frac{1}{2}(\delta t)^{2}f_{i}h_{i,0}-\frac{1}{2}(\delta t)^{2}g_{i}\theta_{i,0}=\frac{1}{2}r^{2}u_{0}(x_{i+1})+(1-r^{2})u_{0}(x_{i})+\frac{1}{2}r^{2}u_{0}(x_{i-1}) \quad \quad \quad \notag \\
&&+(\delta t)v_{0}(x_{i}), \quad \quad  \quad i=\overline{1,(M-1)}.
\label{eqfandgdmsplitefjiszero}
\end{eqnarray}
where $f_{i}:=f(x_{i})$, $h_{i,j}:=h(x_{i},t_{j})$, $g_{i}:=g(x_{i})$ and $\theta_{i,j}:=\theta(x_{i},t_{j})$.

Discretizing \eqref{equofxat0t} and \eqref{equofxatLt}, using \eqref{eq19}, we also have \eqref{equx0} and
\begin{eqnarray}
q_{L}(t_{j})=\frac{\partial{u}}{\partial{x}}(L,t_{j})=\frac{3u_{M,j}-4u_{M-1,j}+u_{M-2,j}}{2(\delta x)},\quad  j=\overline{1,N}. \label{equxL}
\end{eqnarray}

In practice, the additional observations \eqref{equx0} and \eqref{equxL} come from measurement which is inherently contaminated with errors. We therefore model this by replacing the exact data $q_{0}(t)$ and $q_{L}(t)$ by the noisy data \eqref{eq4.5} and
\begin{eqnarray}
q_{L}^{\epsilon}(t_{j})=q_{L}(t_{j})+\tilde{\epsilon}_{j}, \ \ \ j=\overline{1,N},\label{eq4.5uxatLEx5}
\end{eqnarray}
where $(\tilde{\epsilon}_{j})_{j=\overline{1,N}}$ and $N$ random noisy variables generated from a Gaussian normal distribution with mean zero and standard deviation $\tilde{\sigma}=p\times max_{t\in[0,T]}\left|q_{L}(t)\right|$.  

Assembling \eqref{equx0}, \eqref{eqfandgdmsplitef}-\eqref{equxL}, and using \eqref{eq16.1} and \eqref{eq16.2}, the discretised inverse problem reduces to solving a global linear system of $(M-1)\times N+(N+N)$ equations with $(M-1)\times N+((M-1)+(M-1))$ unknowns. Since this system is linear we can eliminate the unknowns $u_{i,j}$ for $i=\overline{1,(M-1)}$, $j=\overline{1,N}$, to reduce the problem to solving an ill-conditioned system of $2N$ equations with $2(M-1)$ unknowns of the form  
\begin{eqnarray}
A(\underline{f},\underline{g})=\underline{b}^{\epsilon}. \label{eqAwithfandg}
\end{eqnarray}
\subsection{Example 5}
This is an example in which we take $c=L=T=1$, $h(x,t)=1$ and $\theta(x,t)=t$ and the input data
\begin{eqnarray}
u(x,0)=u_{0}(x)=\sin(\pi x), \quad
u_{t}(x,0)=v_{0}(x)=x^{2}+1, \quad x\in[0,1],   \label{eq20Ex5}
\end{eqnarray}
\begin{eqnarray}
u(0,t)=P_{0}(t)=t+\frac{t^{2}}{2},\quad u(1,t)=P_{L}(t)=2t+\frac{t^{2}}{2}, \quad t\in(0,1], \label{eq21Ex5}
\end{eqnarray}
\begin{eqnarray}
-\frac{\partial{u}}{\partial{x}}(0,t)=q_{0}(t)=-\pi, \ \ \frac{\partial{u}}{\partial{x}}(1,t)=q_{L}(t)=2t-\pi, \ \ t\in(0,1]. \label{eqq0andqL}
\end{eqnarray}
The exact solution is given by
\begin{eqnarray}
f(x)=1+\pi^{2}\sin(\pi x),\ \ g(x)=-2,\ \ u(x,t)=x^{2}t+\sin(\pi x)+t+\frac{t^{2}}{2},\notag \\ 
(x,t)\in[0,1]\times[0,1].  \label{eq22Ex5}
\end{eqnarray}
  
\begin{figure}[H]
    \begin{subfigure}[b]{0.3\textwidth}
       \caption{}
        \includegraphics[width=11cm]{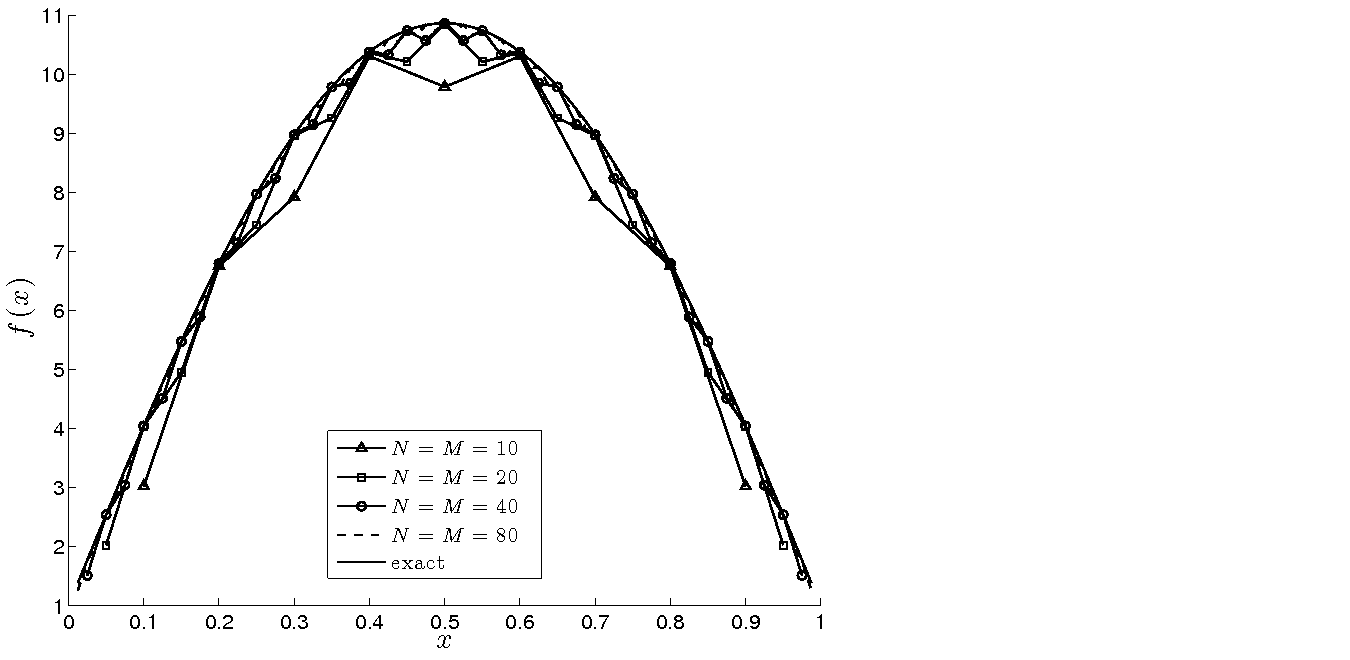}
    \end{subfigure}
    \hspace{3.9cm}
    \begin{subfigure}[b]{0.3\textwidth}
        \caption{}
        \includegraphics[width=11cm]{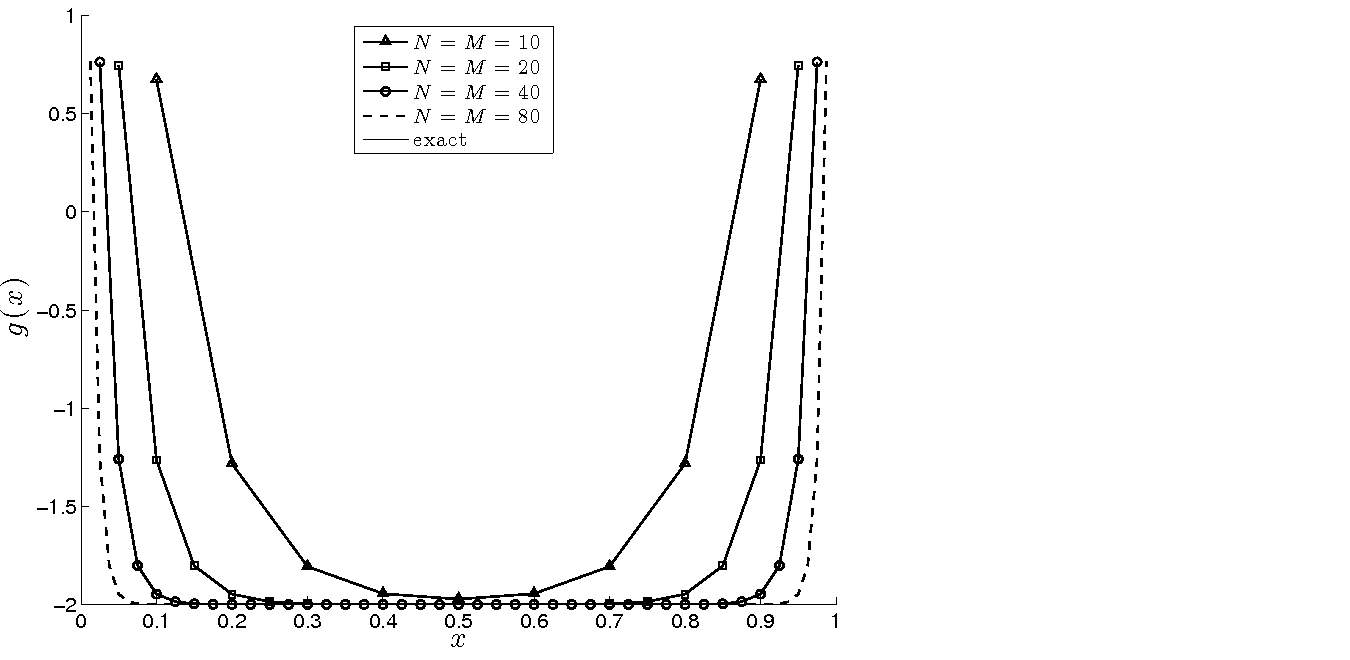}
    \end{subfigure}
\caption{The exact (---) solutions \eqref{eq22Ex5} for the force components $f(x)$ and $g(x)$ in comparison with the numerical solutions for various $N=M\in\lbrace10,20,40,80\rbrace$,  and no regularization, for exact data, for the inverse problem of Example 5.}
\label{fig:inverseproblemexactsolutionoffandgwithnumericalEx5problem3}
\end{figure}

We first consider the case of exact data, i.e. $p=0$ and hence $\underline{\epsilon}=\underline{\tilde{\epsilon}}=\underline{0}$ in \eqref{eq4.5} and \eqref{eq4.5uxatLEx5}. The numerical results corresponding to $f(x)$ and $g(x)$ are plotted in Figure \ref{fig:inverseproblemexactsolutionoffandgwithnumericalEx5problem3}. From this figure it can be seen that convergent and accurate numerical solutions are obtained especially, for $f(x)$, although for $g(x)$ are some inaccuracies manifested near the end points $x\in\lbrace0,1\rbrace$.

\begin{figure}[H]
    \begin{subfigure}[b]{0.3\textwidth}
       \caption{}
        \includegraphics[width=11cm]{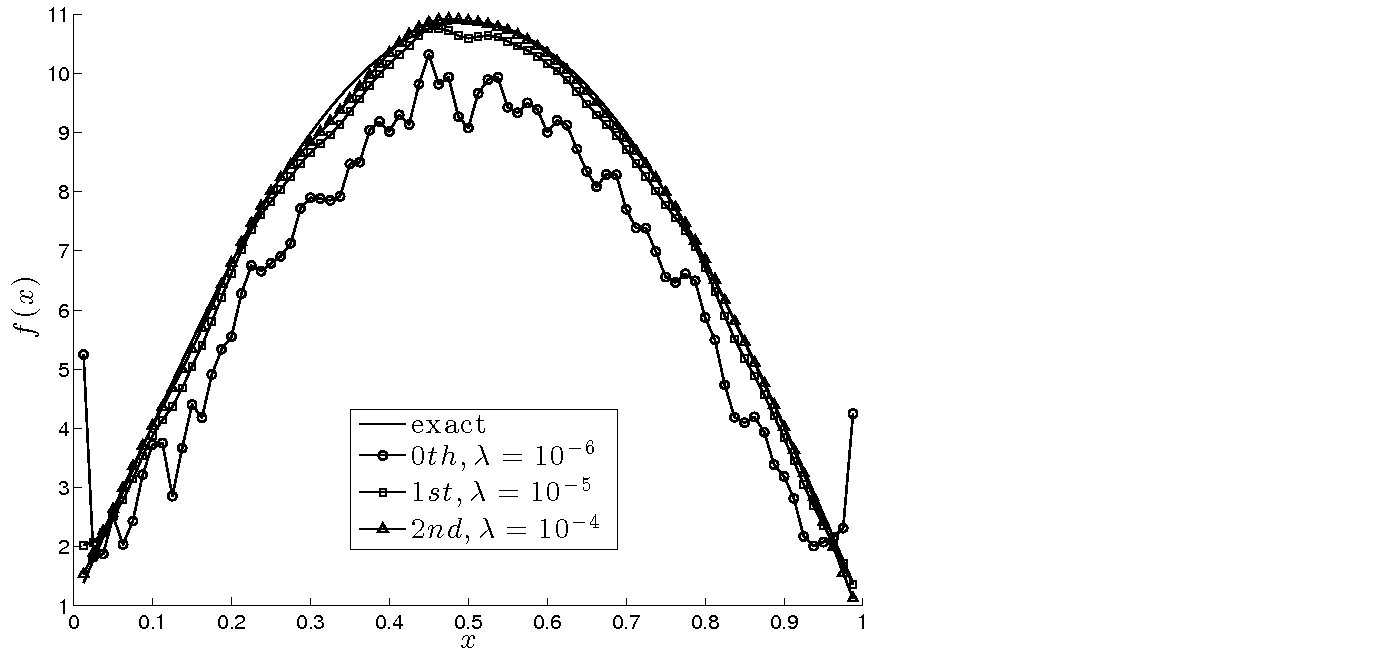}
    \end{subfigure}
    \hspace{3.9cm}
    \begin{subfigure}[b]{0.3\textwidth}
        \caption{}
        \includegraphics[width=11cm]{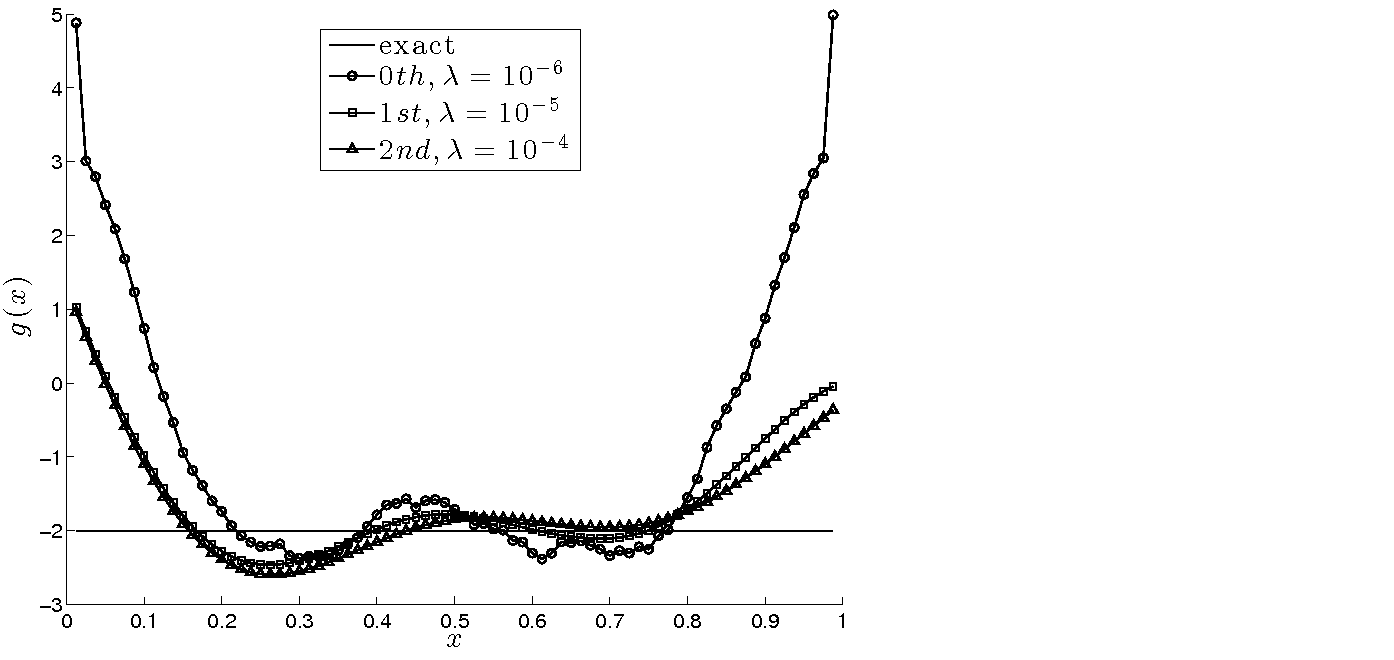}
    \end{subfigure}
\caption{The exact (---) solutions \eqref{eq22Ex5} for the force components $f(x)$ and $g(x)$ in comparison with the numerical solutions, for $N=M=80$, $p=1\%$ noise and various order regularization methods, for the inverse problem of Example 5.}
\label{fig:optimalfandgwithexactin0th1st2ndEx5problem3}
\end{figure}

We include some ($p=1\%$) noise into the input data \eqref{equx0} and \eqref{equxL}, as given by equations \eqref{eq4.5} and \eqref{eq4.5uxatLEx5}. Figure \ref{fig:optimalfandgwithexactin0th1st2ndEx5problem3} shows the regularized numerical solutions for $f(x)$ and $g(x)$ obtained with various regularizations and one can observe that reasonably stable numerical solutions are obtained. 
\section{Conclusions}
In this paper, the determination of space-dependent forces from boundary Cauchy data in the wave equation has been investigated. The solution of this linear inverse problem is unique, but is still ill-posed since small errors in the input flux cause large errors in the output force. The problem is discretised numerically using the FDM, and in order to stabilise the solution, the Tikhonov regularization method has been employed. The choice of the regularization parameter was based on the L-curve criterion. 
Numerical examples indicate that the method can accurately recover the unknown space-dependent force. The time-dependent force identification will be investigated in Part II, \cite{husseinlesnic}.\\
\\  
\textbf{\large Acknowledgments}\\
S.O. Hussein would like to thank the Human Capacity Development Programme (HCDP) in Kurdistan for their financial support in this research.

\end{document}